\documentclass[11pt]{article}

\usepackage[margin=1in]{geometry}

\usepackage{algorithm, algpseudocode}
\usepackage{bm,amsmath,commath,mathtools,amsfonts,amssymb}
\usepackage{graphicx, caption, subcaption, float}
\usepackage{url, cite, hyperref}
\hypersetup{
    colorlinks=true,
    linkcolor=blue,
    citecolor = red,
}

\numberwithin{equation}{section}
\usepackage{authblk}

\usepackage{tikz}
\usetikzlibrary{arrows,positioning}
\usepackage{onimage}

\tikzstyle{block} = [draw, rectangle, align=center, minimum width=2cm, inner sep=1ex]

\usepackage{color}

\title{\textbf{Evaluating the accuracy of the dynamic mode decomposition}}

\author[1]{Hao Zhang\footnote{Email: haozhang@princeton.edu}}
\author[1]{Scott T. M. Dawson}
\author[1]{Clarence W. Rowley}
\author[2]{Eric A. Deem}
\author[2]{Louis N. Cattafesta}
\affil[1]{Mechanical and Aerospace Engineering, Princeton University}
\affil[2]{Mechanical Engineering, Florida State University}

\begin{document}

\maketitle

\begin{abstract}
Dynamic mode decomposition (DMD) gives a practical means of extracting dynamic information from data, in the form of spatial modes and their associated frequencies and growth/decay rates.
DMD can be considered as a numerical approximation to the Koopman operator, an infinite-dimensional linear operator defined for (nonlinear) dynamical systems.
This work proposes a new criterion to estimate the accuracy of DMD on a
mode-by-mode basis, by estimating how closely each individual DMD eigenfunction
approximates the corresponding Koopman eigenfunction.  This approach does not require any prior knowledge of the system dynamics or the true Koopman spectral decomposition.
The method may be applied to extensions of DMD (i.e., extended/kernel DMD), which are applicable to a wider range of problems.
The accuracy criterion is first validated against the true error with a
synthetic system for which the true Koopman spectral decomposition is known.
We next demonstrate how this proposed accuracy criterion can be used to assess the performance of various choices of kernel when using the kernel method for extended DMD.
Finally, we show that our proposed method successfully identifies modes of high
accuracy when applying DMD to data from experiments in fluids, in particular
particle image velocimetry of a cylinder wake and a canonical separated boundary
layer.

\end{abstract}

\section{Introduction}

The decomposition of spatiotemporal data into spatial modes and temporal functions describing their evolution gives a means to isolate coherent features and assemble low-order representations of complex dynamics.
Over the past decade, the dynamic mode decomposition (DMD) \cite{Schmid2010DMD} has become a routinely-used method for such purposes \cite{rowley2009spectral,tu2014spectral,DeemAIAA2017,zhang2017online}. See, for example, \cite{kutz2016book} and \cite{rowley2017arfm} for reviews of many ensuing uses and applications of DMD.
 While successfully used on a range of datasets, general questions still exist in terms of how to select a reduced set of modes, and how to ensure results are quantitatively accurate.
 On the first point, numerous methods have been proposed to select a reduced number of modes that best represent the dynamics of the system  \cite{chen2012variants,wynn2013optimal,jovanovic2014sparsity,kou2016improved}.
On the second point, the sensitivity of the outputs of DMD to noisy data has also been investigated \cite{duke2012error,Bagheri2014noise}, and a number of modified algorithms proposed that give improved accuracy for noisy data \cite{hemati2015biasing,dawson2016characterizing}.

The present work differs from these past studies by giving a means of estimating the accuracy of DMD on a mode-by-mode basis, without any a-priori knowledge of the system dynamics, noise characteristics, or truncation of low-energy modes.
It has been shown previously \cite{rowley2009spectral,tu2014dynamic} that DMD approximates the Koopman operator, an infinite-dimensional linear operator defined for (nonlinear) dynamical systems. In this work, we will exploit this connection by estimating the accuracy to which we approximate eigenfunctions of the Koopman operator.
This approach allows our analysis to naturally extend to extensions of DMD \cite{williams2015data} that are designed to improve the approximation to the Koopman operator for nonlinear systems. Extended DMD uses nonlinear observables to expand the space in which the Koopman operator is approximated. However, EDMD suffers from the curse of dimensionality: that is, the computational cost increases rapidly with the dimension of the state. To circumvent this issue, kernel DMD (KDMD) \cite{williams2014kernel} was proposed as a computationally inexpensive alternative, which makes use of a kernel function to implicitly include a rich (and nonlinear) set of observables, while maintaining the same computational cost as DMD. The optimal choice of kernel function for KDMD is still an open question, and here we demonstrate that the accuracy criterion may be used to evaluate and compare the performance of various kernels.

The structure of this work is as follows. We first review DMD, the Koopman operator, and kernel DMD in section~\ref{sec:background}, before presenting and validating our proposed accuracy criterion in section~\ref{sec:crit}. Section~\ref{sec:synth} uses the accuracy criterion to measure the performance of various kernels in KDMD for a simple nonlinear system, while section \ref{sec:exp} demonstrates that this criterion is effective in selecting accurate DMD modes from experimental data.

\section{Background}
\label{sec:background}
We first give a review of previous results, including the DMD algorithm and its
connections to the Koopman operator (section~\ref{sec:DMD}), as well as extensions of
DMD that can better approximate the Koopman operator for nonlinear systems (section~\ref{sec:kdmd}).

\subsection{Dynamic mode decomposition}
\label{sec:DMD}
Dynamic mode decomposition was introduced in~\cite{Schmid2010DMD}, and our
presentation here follows that in~\cite{tu2014dynamic,rowley2017arfm}.  Consider a discrete-time dynamical system whose state space is denoted by
$X\subset \mathbb{R}^n$, and suppose the dynamics are given by
\begin{equation}
\label{eq:dynamics}
\bm{x}({k+1}) = \bm{F}(\bm{x}({k})), \quad \bm{x}(k) \in X.
\end{equation}
Let $\psi_1,\ldots,\psi_q$ be real-valued functions on $X$, which we call
observables, and let $\bm{\psi}: X \to \mathbb{R}^q$ denote the
vector-valued function whose components are $(\psi_1,\ldots,\psi_q)$. We may not
be able to measure the state $\bm{x}$ directly, but instead, we can measure the
vector
\begin{equation*}
\bm{y} = \bm{\psi}(\bm{x}).
\end{equation*}
As a special case, $\bm{y}$ could be the state itself, i.e., $\bm{y} = \bm{\psi}(\bm{x}) = \bm{x}$. For complex systems, it can be advantageous to define observables that are nonlinear functions of the state, which will be discussed in more detail in section~\ref{sec:kdmd}. For the purposes of describing standard DMD, we assume $\bm{y} = \bm{x}$.

We consider pairs of snapshots $(\bm{x}_k,\bm{x}^\#_k)$, with $ \bm{x}_k \in X,
k = 1,2,\cdots,m$, and where $\bm{x}^\#_k=F(\bm{x}_k)$ is the image of
$\bm{x}_k$ upon application of the dynamics  \eqref{eq:dynamics}. For sequential
data, $\bm{x}(1),\ldots,\bm{x}(m+1)$ satisfying~\eqref{eq:dynamics}, one takes
$\bm{x}_k=\bm{x}(k)$, $\bm{x}^\#_k=\bm{x}({k+1})$, though non-sequential data
may also be used, such as from multiple runs of experiments or simulations
\cite{tu2014dynamic}. In DMD, we seek a matrix $\bm{A} \in \mathbb{R}^{q \times q}$ such that
\begin{equation*}
\bm{y}^\#_k = \bm{A} \bm{y}_k, \ \  k=1,2,\cdots,m
\end{equation*}
holds, at least approximately. We form two matrices
\begin{equation*}
\bm{Y} =
\begin{bmatrix}
\bm{y}_1 & \bm{y}_2 & \cdots & \bm{y}_m
\end{bmatrix}, \quad
\bm{Y}^\# =
\begin{bmatrix}
\bm{y}^\#_1 & \bm{y}^\#_2 & \cdots & \bm{y}^\#_m
\end{bmatrix},
\end{equation*}
and define the DMD matrix $\bm{A}$ by
\begin{equation}
\label{eq:DMDmatrix}
\bm{A} = \bm{Y}^\# \bm{Y}^+.
\end{equation}
DMD modes and eigenvalues are the eigenvectors and eigenvalues of  $\bm{A}$. A typical algorithm to compute these modes and eigenvalues is as follows\cite{tu2014dynamic}:
\newline

\textbf{Algorithm (DMD)}
\begin{enumerate}
\item Compute the reduced SVD $\bm{Y}=\bm{U}\bm{\Sigma}\bm{V}^T$.
\item (Optional) Truncate the SVD by only retaining the first $r$ colunms of $\bm{U},\bm{V}$, and the first $r$ rows and columns of $\bm{\Sigma}$, to obtain $\bm{U}_r,\bm{\Sigma}_r,\bm{V}_r$.
\item Let $\tilde{\bm{A}}=\bm{U}_r^T \bm{A} \bm{U}_r=\bm{U}_r^T \bm{Y}^\# \bm{V}_r \bm{\Sigma}_r^{-1}$, $\tilde{\bm{A}} \in \mathbb{R}^{r \times r}$.
\item Find the eigenvalues $\mu_i$ and eigenvectors $\tilde{\bm{v}}_i$ of $\tilde{\bm{A}}$, such that $\tilde{\bm{A}} \tilde{\bm{v}}_i = \mu_i \tilde{\bm{v}}_i$.
\item The (projected) DMD modes are given by $\bm{v}_i=\bm{U}_r^T \tilde{\bm{v}}_i$, with corresponding (discrete-time) DMD eigenvalues $\mu_i$.
\end{enumerate}

The eigenvectors of the matrix $\bm{A} \in \mathbb{R}^{q \times q}$ can be
found from the eigenvectors of the smaller
matrix~$\tilde{\bm{A}}\in\mathbb{R}^{r\times r}$. We denote the eigenvalues and
eigenvectors of $\bm{A}$ by $\{\mu_i,\bm{v}_i\}$. In the case of sequential data
(for which $\bm{y}^\#_k=\bm{y}_{k+1}$), suppose that we can express the initial state as
\begin{equation*}
\bm{y}_{1} = \sum_{i=1}^q c_i \bm{v}_i.
\end{equation*}
The time evolution of the system (starting at $\bm{y}_1$) is then predicted by DMD to be
\begin{equation}
\label{eq:DMDevolution}
\bm{y}_{k+1} = \bm{A}^k \bm{y}_1 = \sum_{i=1}^q c_i \mu_i^k \bm{v}_i.
\end{equation}
Therefore, each DMD mode $\bm{v}_i$ is associated with a single frequency and growth/decay rate (DMD eigenvalue $\mu_i$). In reality, \eqref{eq:DMDevolution} may not hold exactly, depending on the quantity and quality of data used, whether the system dynamics are nonlinear, whether the SVD is truncated in step 2 of the DMD algorithm above.
 For cases where equation \eqref{eq:DMDevolution} does not give an exact description of the dynamics, DMD gives a least-squares fit to the data (as pairs of snapshots).

There are connections between DMD and an infinite-dimensional linear operator called the Koopman
operator \cite{tu2014dynamic, rowley2009spectral}, with the high-level idea
being that DMD gives a finite-dimensional numerical approximation of the Koopman
operator. Our proposed criterion for evaluating the accuracy of DMD
exploits this connection.  For a given state-space~$X$, the Koopman operator
acts on scalar-valued functions of~$X$, which we referred earlier as
observables.  Here, we consider observables in $L^2(X)$, the space of square integrable
functions on~$X$.  Given the dynamics in~\eqref{eq:dynamics}, one then defines the Koopman operator\footnote{To be fully rigorous, one typically assumes the dynamics~\eqref{eq:dynamics}
are measure-preserving, so that $\phi\circ \bm{F}$ is in $L^2$ whenever $\phi\in
L^2$; in fact, if $\mathbf{F}$ is measure preserving, then $\mathcal{K}$ is an isometry.}
 $\mathcal{K}:L^2(X)\to L^2(X)$ by
\begin{equation}
(\mathcal{K} \phi)(\bm{x}) = (\phi \circ \bm{F}) (\bm{x}) = \phi(\bm{F}(\bm{x})).
\end{equation}
That is, $\mathcal{K}$ maps a function $\phi$ to another function
$\phi\circ\bm{F}$, and $(\mathcal{K}\phi)(\bm{x})$ gives the value of $\phi(\bm{x})$ at the next time step.  Here we emphasize two points: first that the Koopman operator acts on functions of the state instead of the state itself; and second that the Koopman operator is linear, even though the dynamics might be nonlinear. On the second point, note that
\begin{equation*}
\mathcal{K}(c_1 \phi_1 + c_2 \phi_2)=c_1 \mathcal{K} \phi_1 + c_2 \mathcal{K} \phi_2
\end{equation*}
holds for any functions $\phi_1, \phi_2$ and any scalars $c_1, c_2$. Since the Koopman operator is linear, it may have eigenvalues and eigenfunctions, which satisfy
\begin{equation}
\mathcal{K} \varphi=\mu \varphi,
\end{equation}
where $\varphi$ is the eigenfunction with eigenvalue $\mu$.



Now, suppose we have a given set of observables $\{ \psi_1, \psi_2, \cdots,
\psi_q \}$, and suppose $\varphi$ is a Koopman eigenfuntion (with eigenvalue $\mu$) that lies in the span of $\{\psi_j\}$: i.e.,
\begin{equation}
\varphi(\bm{x})=\bar{w}_1\psi_1(\bm{x})+\cdots+\bar{w}_q\psi_q(\bm{x})=\bm{w}^* \bm{\psi}(\bm{x}),
\end{equation}
for some $\bm{w}^* \in \mathbb{C}^n$. Then one can show (see \cite[\S
4.1]{tu2014dynamic}) that under certain conditions on the data, $\bm{w}^*$ is a left eigenvector of the DMD matrix $\bm{A}$ with eigenvalue $\mu$ (i.e., $\bm{w}^*\bm{A}=\mu \bm{w}^*$). This connection implies that we can approximate Koopman eigenfunctions (and eigenvalues) for a given unknown dynamical system directly from data using DMD. In particular, given left eigenvectors of the DMD matrix ($\bm{w}_i^*\bm{A}=\mu_i \bm{w}_i^*$), we consider $\varphi_i(\bm{x})=\bm{w}_i^* \bm{\psi}(\bm{x})$ as a DMD-approximated Koopman eigenfunction, with eigenvalue $\mu_i$.

\subsection{Extended DMD and kernel DMD}
\label{sec:kdmd}
In order to apply the connection between DMD and Koopman mentioned above, the
Koopman eigenfunctions must lie within the space spanned by the observables
$\{\psi_j\}$. If one takes $\bm{\psi}(\bm{x})=\bm{x}$, as with standard DMD, then
the subspace spanned by $\{\psi_j\}$ consists only of linear functions
of~$\bm{x}$, and this subspace is often not large enough to include
eigenfunctions of~$\mathcal{K}$ (a notable exception being the case in which
$\bm{F}$ is linear).  Extended DMD (EDMD) was proposed
in~\cite{williams2015data} in order to enlarge
the subspace of observables, and therefore better approximate Koopman
eigenfunctions. In particular, Extended DMD approximates the Koopman operator by
a weighted residual method, with trial functions given by $\{\psi_j\}$ and a
particular choice of test
functions specified by the data.  Examples of observables $\psi_j(\bm{x})$ could include
polynomials, Fourier modes, indicator functions, or spectral elements, as
suggested in \cite{williams2015data}. For instance, if we take $\bm{x} \in
\mathbb{R}^2$ and take observables to be monomials in components of $\bm{x}$ up to degree $d=2$ (including the constant 1), then the vector of observables is
\begin{equation*}
\bm{\psi}(\bm{x})=
\begin{bmatrix}
1 & x_1 & x_2 & x_1^2 & x_1 x_2 & x_2^2
\end{bmatrix}^T.
\end{equation*}
We can potentially approximate many more accurate Koopman eigenfunctions with
EDMD than we could with DMD. However, EDMD suffers from the curse of
dimensionality \cite{bishop2006pattern}. If the state dimension is $n$
and we consider (multivariate) polynomials up to degree~$d$, then the number of
observables is $q=\binom{n+d}{d}$, which is approximately $n^d$ for large~$n$.  For large
problems (as arise in fluids), data might typically have $n \approx 10^6$, so
even if one considers only quadratic polynomials, the number of observables is $q\approx
10^{12}$, too large for practical computation. It is thus very computationally
expensive to consider large subspaces of observables.

Kernel DMD (KDMD) has been proposed to deal with this curse of dimensionality \cite{williams2014kernel}. In KDMD, EDMD is reformulated such that only inner products of observables need to be computed. The inner product can be evaluated by making use of a kernel function, a common technique in the community of machine learning. A kernel function $k: \mathbb{R}^n \times \mathbb{R}^n \rightarrow \mathbb{R}$ is defined as
\begin{equation}
k(\bm{x},\hat{\bm{x}})=\langle \bm{\psi}(\bm{x}), \bm{\psi}(\hat{\bm{x}}) \rangle.
\end{equation}
To appreciate how kernel functions work, consider for example a polynomial kernel $k(\bm{x},\hat{\bm{x}})=(1+\bm{x}^T\hat{\bm{x}})^d$ as an example. This kernel corresponds to a set of observables $\bm{\psi}(\bm{x})$ consisting of all monomials in components of $\bm{x}$ up to degree $d$ \cite{bishop2006pattern}. Taking $n=2$ and $d=2$, this kernel function can be expanded as
\begin{equation} \label{eq:polykerneld2}
\begin{aligned}
(1+\bm{x}^T\hat{\bm{x}})^2 & =1+2x_1 \hat{x}_1 + 2x_2 \hat{x}_2 + 2x_1^2 \hat{x}_1^2 + 2 x_1 x_2 \hat{x}_1 \hat{x}_2 + \hat{x}_1^2 \hat{x}_2^2 \\
 & =  \langle \bm{\psi}(\bm{x}),\bm{\psi}(\hat{\bm{x}}) \rangle,
\end{aligned}
\end{equation}
where
$\bm{\psi}(\bm{x})=(1,\sqrt{2}x_1,\sqrt{2}x_2,x_1^2,\sqrt{2}x_1x_2,x_2^2)$. In
the terminology of machine learning, $\bm{\psi}$ is called the feature map, and
$\bm{\psi}(\bm{x}) \in \mathbb{R}^q$ is called the feature space (which might be
infinite dimensional). In the example above, the dimension of the (implicitly
defined) feature space is $q=6$, but in order to compute $k(\bm{x},\hat{\bm{x}})$,
we require inner products only in state space, which has dimension $n=2$. Kernel functions hence can be used to evaluate the inner product in a high dimensional (or even infinite dimensional) feature space in an efficient way. More examples of kernel functions are given in section~\ref{sec:kernelfunc}.

\section{Accuracy criterion for DMD}
\label{sec:crit}

The connection between DMD and the Koopman operator as discussed in
section~\ref{sec:DMD} implies that we can use variants of DMD (e.g., DMD, EDMD,
or KDMD) to approximate Koopman eigenfunctions and
eigenvalues, given access to
data. By applying DMD variants to a given dataset, we can potentially identify
many Koopman eigenfunctions and eigenvalues (which we refer to as
eigenpairs). However, the reliability of these eigenpairs remains
unknown. Before using DMD results for any analysis or reduced order modeling, it
is desirable and necessary to assess the quality (i.e., accuracy) of the
results. In this section, we will develop a criterion for evaluating the
accuracy of DMD-approximated Koopman eigenpairs. We describe this accuracy
criterion in section \ref{sec:proposemetric}, and then validate it in section
\ref{sec:Validating} using a simple nonlinear system where the analytical
Koopman eigenpairs are known.

The most common way to select which of the computed DMD modes are most relevant
is to use the ``mode amplitude'': for sequential data, one projects the initial
condition onto DMD modes and one views the magnitude of the projection
coefficients as the mode amplitudes.  It is common practice
\cite{tu2014dynamic,hemati2016improving,rowley2009spectral} to retain the modes
of largest amplitude.  This approach sounds plausible; however, it was observed
in \cite{jovanovic2014sparsity} (which used sparsity-promoting techniques to
select modes) that mode amplitude is not always a useful criterion for mode
selection.  Indeed, mode amplitudes can be misleading, as we illustrate below with a
simple example.

Suppose we have three DMD modes,
\begin{equation}
\bm{v}_1 = (1,0,0),\qquad
\bm{v}_2 = (0,1,0),\qquad
\bm{v}_3 = (0,1,\epsilon),
\end{equation}
where $\epsilon$ is small and thus $\bm{v}_2$ and $\bm{v}_3$ are almost parallel. If we consider an initial condition $\bm{x}_0= (1,0,\zeta)$, and project it onto these DMD modes, we obtain
\begin{equation}
\bm{x}_0 = \bm{v}_1 -\frac{\zeta}{\epsilon} \bm{v}_2 + \frac{\zeta}{\epsilon} \bm{v}_3.
\end{equation}
For instance, if $\zeta=10^{-3}$ and $\epsilon=10^{-6}$, then $\zeta/\epsilon =
10^3$, so the mode amplitude (defined as the magnitude of the projection
coefficients) indicates that $\bm{v}_2$ and $\bm{v}_3$ are much more important
than $\bm{v}_1$. The mode amplitudes indicate that we might be able to neglect
$\bm{v}_1$ without significant adverse effects. However, it is clear that
$\bm{v}_1$ is much more relevant for reconstructing $\bm{x}_0$: if we use only
$\bm{v}_2$ and~$\bm{v}_3$, we obtain
\begin{equation}
-\frac{\zeta}{\epsilon} \bm{v}_2 + \frac{\zeta}{\epsilon} \bm{v}_3 = (0,0,\zeta),
\end{equation}
which does not accurately approximate~$\bm{x}_0=(1,0,\zeta)$.  A better
approximation to $\bm{x}_0$ is simply $\bm{v}_1=(1,0,0)$.  This example
illustrates that mode amplitude is not always a reliable criterion for selecting
which modes are important, especially when modes are almost parallel.  Note that
this problem would also arise if using other methods to measure mode amplitude
(e.g., \cite{kou2016improved}).

The accuracy criterion we describe below does not provide a way of selecting
which modes are dominant, and in the example above, accuracy of the modes did
not play a role.  However, the accuracy criterion can provide a way of
eliminating candidate modes that we know to be inaccurate, so in this way it can
help with the problem of mode selection.

\subsection{Proposed accuracy criterion}
\label{sec:proposemetric}

Given data from an experiment or simulation, we can split the dataset into
training data and testing data. Training data is used to approximate DMD modes
(and associated Koopman eigenpairs), while testing data is used to evaluate the
quality of these identified modes. Data-driven algorithms may suffer from the
problem of over-fitting \cite{hawkins2004problem}, so any evaluation criteria
should use testing data that differs from the training data.

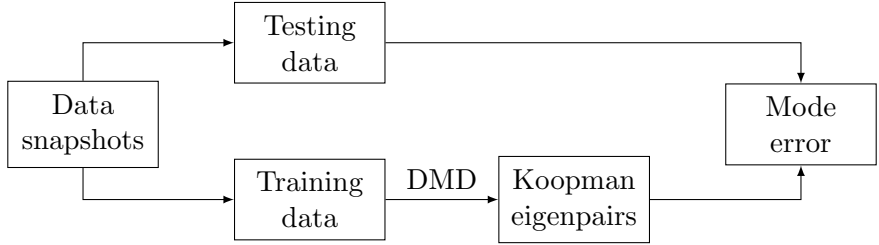
\begin{figure}
  \centering
  \begin{tikzpicture}[node distance=1cm and 1.5cm,>=latex]
    \node[block] (train) {Training\\data};
    \node[block, above=of train] (test) {Testing\\data};
    \node[block, left=1cm of train, yshift=1cm] (data) {Data\\snapshots};
    \node[block, right=of train] (koopman) {Koopman\\eigenpairs};
    \node[block, right=1cm of koopman, yshift=1cm] (error) {Mode\\error};
    \draw[->] (data) |- (test);
    \draw[->] (data) |- (train);
    \draw[->] (train) -- node[above] {DMD} (koopman);
    \draw[->] (test) -| (error);
    \draw[->] (koopman) -| (error);
  \end{tikzpicture}
  \caption{A diagram summarizing the implementation of the accuracy
    criterion. Training data is used to approximate Koopman eigenpairs with
    variants of DMD, while testing data is used to evaluate the quality of Koopman
    eigenpairs.}
  \label{fig:accuracy}
\end{figure}

The idea of our approach is to evaluate the accuracy of a DMD mode (and
eigenvalue) by looking at the accuracy of its corresponding Koopman
eigenfunction.  Suppose we are given an approximate Koopman eigenpair
$(\mu,\varphi)$, and we wish to evaluate its accuracy.  If $(\mu,\varphi)$ were
a true Koopman eigenpair, then by definition it would satisfy
\begin{equation*}
\varphi\circ \bm{F} = \mu \varphi,
\end{equation*}
where $\bm{F}$ defines the dynamics in~\eqref{eq:dynamics}. Ideally, we would like to compute
\begin{equation}
\label{eq:errorfuncnorm}
\frac{\|\varphi\circ \bm{F}- \mu \varphi\|}{\|\varphi\|},
\end{equation}
where $\|\cdot\|$ is the norm of a function. (We divide by $\|\varphi\|$ so that
the above quantity is independent of the scaling of the
eigenfunction~$\|\varphi\|$.)  However, in order to
compute~\eqref{eq:errorfuncnorm}, we require
explicit knowledge of the dynamics $\bm{F}$, which is unknown in most
cases of interest. Instead, we can estimate the above quantity using finite
number of data points (i.e., the testing data). The estimation should give some
sense of the quantity in~\eqref{eq:errorfuncnorm}, using only the testing data,
which consists of pairs of samples $(\bm{x}_k,\bm{x}_k^\#)$ with
$\bm{x}_k\in X$ and $\bm{x}_k^\#=\bm{F}(\bm{x}_k)$.  This observation motivates the following definition of an accuracy criterion:
\begin{equation} \label{eq:error}
\alpha = \frac{\sum_k |\varphi(\bm{x}_k^\#) - \mu \varphi(\bm{x}_k)|}{\sum_k |\varphi(\bm{x}_k)|},
\end{equation}
where $|\cdot|$ denotes the absolute value, and the summation is over the entire testing dataset. A diagram summarizing how this accuracy criterion may be applied is shown in Figure \ref{fig:accuracy}. More specifically, given a DMD-approximated eigenfunction $\varphi(\bm{x})= \bm{w}^* \bm{\psi}(\bm{x})$ with eigenvalue $\mu$ (i.e.,  $\bm{w}^*\bm{A}=\mu \bm{w}^*$, with $\bm{A}$ as defined in~\eqref{eq:DMDmatrix}), the accuracy criterion, or estimated mode error, can be written as
\begin{equation}
\alpha = \frac{\sum_k |\bm{w}^* \bm{\psi}(\bm{x}_k^{\#}) - \mu \bm{w}^* \bm{\psi}(\bm{x}_k)|}{\sum_k |\bm{w}^* \bm{\psi}(\bm{x}_k)|}.
\end{equation}

The numerator measures to what extent the eigenfunction equation holds, and the denominator gives a measure of the magnitude of the eigenfunction. Here $\alpha$ can be interpreted as the error of a Koopman eigenpair. The error is defined on a mode-by-mode basis, which enables independent evaluation for each individual DMD mode. Therefore it makes sense to call $\alpha$ the mode error. Observe that $\alpha$ is always non-negative, and it is usually less than 1. When we feed in the true Koopman eigenfunction and eigenvalue into $\alpha$ in equation (\ref{eq:error}), then $\alpha=0$ (assuming that the testing data is noise-free).
If $\alpha$ is close to~1, the Koopman eigenpair is extremely unreliable,
because the discrepancy in the eigenfunction equation is of the same order as
the magnitude of the eigenfunction. Therefore, usually we only care about the
DMD eigenpairs for which $0\le\alpha\ll 1$. In our definition in \eqref{eq:error}, we have used the absolute value to indicate the discrepancy in the eigenfunction equation. However, it is also possible to use other norms, such as the $\ell^2$ norm (or its square), which yield similar results in terms of indicating relative accuracy of modes.

A meaningful evaluation criterion should be (fairly) independent the scaling of
the eigenfunctions, the scaling of the testing data, and the size of the testing
set.
The proposed accuracy criterion approximately satisfies all of these. To show this, we consider the simple case where the full system state is used in DMD, i.e., $\bm{\psi}(\bm{x})=\bm{x}$ and the DMD-computed eigenfunction is linear, i.e., $\varphi(\bm{x}) = \bm{w}^*\bm{\psi}(\bm{x})=\bm{w}^*\bm{x}$.
The fact that we normalize by the magnitude of the observables means that $\alpha$ is relatively independent of eigenfunction scaling, data scaling, and data quantity, as is desired.
In the case where the observable is not the full state (i.e., when using EDMD or KDMD), the scaling of the eigenfunctions and size of testing again do not influence $\alpha$, for the same reason. However, due to the nonlinear transformation $\bm{\psi}(\bm{x})$, the scaling of testing data $\bm{x}$ may play some role in the size of $\alpha$. Fortunately, it is reasonable to expect that the relative magnitude of $\alpha$ should still indicate the relative accuracy of different DMD-computed Koopman eigenpairs.


We point out that if the testing data is clean, mode error is determined only
by the quality of DMD-approximated Koopman eigenpairs. If the testing data is
noisy, mode error is also affected by the noise in testing data. For
experimental data, we have access only to the noisy measurements. In these
cases, the relative magnitude of~$\alpha$ is still expected to indicate relative
accuracy of Koopman eigenpairs. We reiterate again that this definition of error
does not assume access to analytical Koopman spectral decomposition, which is
unknown in most cases.


\subsection{Validating the accuracy criterion}
\label{sec:Validating}
We have proposed an accuracy criterion that exploits the connection between DMD
and the Koopman operator. Before applying this criterion to real data, we first
seek to validate it as a reliable measure of accuracy. We will first consider a
simple 2D nonlinear system for which the analytical Koopman spectral
decomposition is known. Given analytical Koopman eigenpairs, we can define the
true error to be the distance between the DMD eigenvalue and the true eigenvalue
(eigenvalue error), or the difference between the DMD eigenfunction and the true eigenfunction (eigenfunction error). We will validate the accuracy criterion against the true error, and show that accuracy criterion reliably indicates accuracy.

Here we consider a 2D nonlinear map (also considered in \cite{tu2014dynamic}) with dynamics defined by
\begin{equation}
\label{eq:2Dex}
\begin{bmatrix}
x_1 \\
x_2
\end{bmatrix}
\mapsto
\begin{bmatrix}
\gamma x_1 \\
\delta x_2 + (\gamma^2-\delta)x_1^2
\end{bmatrix},
\quad \gamma = 0.9, \delta = 0.8.
\end{equation}
%
It is straightforward to verify that $\gamma,\delta$ are Koopman eigenvalues with respective eigenfunctions
\begin{equation*}
\varphi_{\gamma}(\bm{x}) = x_1,\qquad
\varphi_{\delta}(\bm{x}) = x_2 - x_1^2.
\end{equation*}
Additional Koopman eigenvalues and eigenfunctions are given by
\begin{equation} \label{eq:2Danal}
{\mu}_{k,\ell}  =\gamma^k \delta^{\ell}, \qquad {\varphi}_{k,\ell}  =\varphi^k_{\gamma} \varphi^{\ell}_{\delta},
\end{equation}
where $k, \ell=0,1,2,\cdots$ are non-negative integers. Notice that the
analytical eigenfunctions are multivariate polynomials in the state variables.

To collect training data, $m=100$ random initial points are sampled from a uniform distribution on $[-1,1] \times [-1,1]$, and their images are found by applying the map  defined in equation~\eqref{eq:2Dex}. Similarly we also generate $m_{\text{test}}=100$ snapshot pairs as the testing data. The generated training and testing dataset are used for subsequent analysis in both this and the next section.

\begin{figure}[t]
        \centering
        \begin{subfigure}[b]{0.45\textwidth}
    \centering
    \begin{tikzonimage}[width=.9\linewidth]{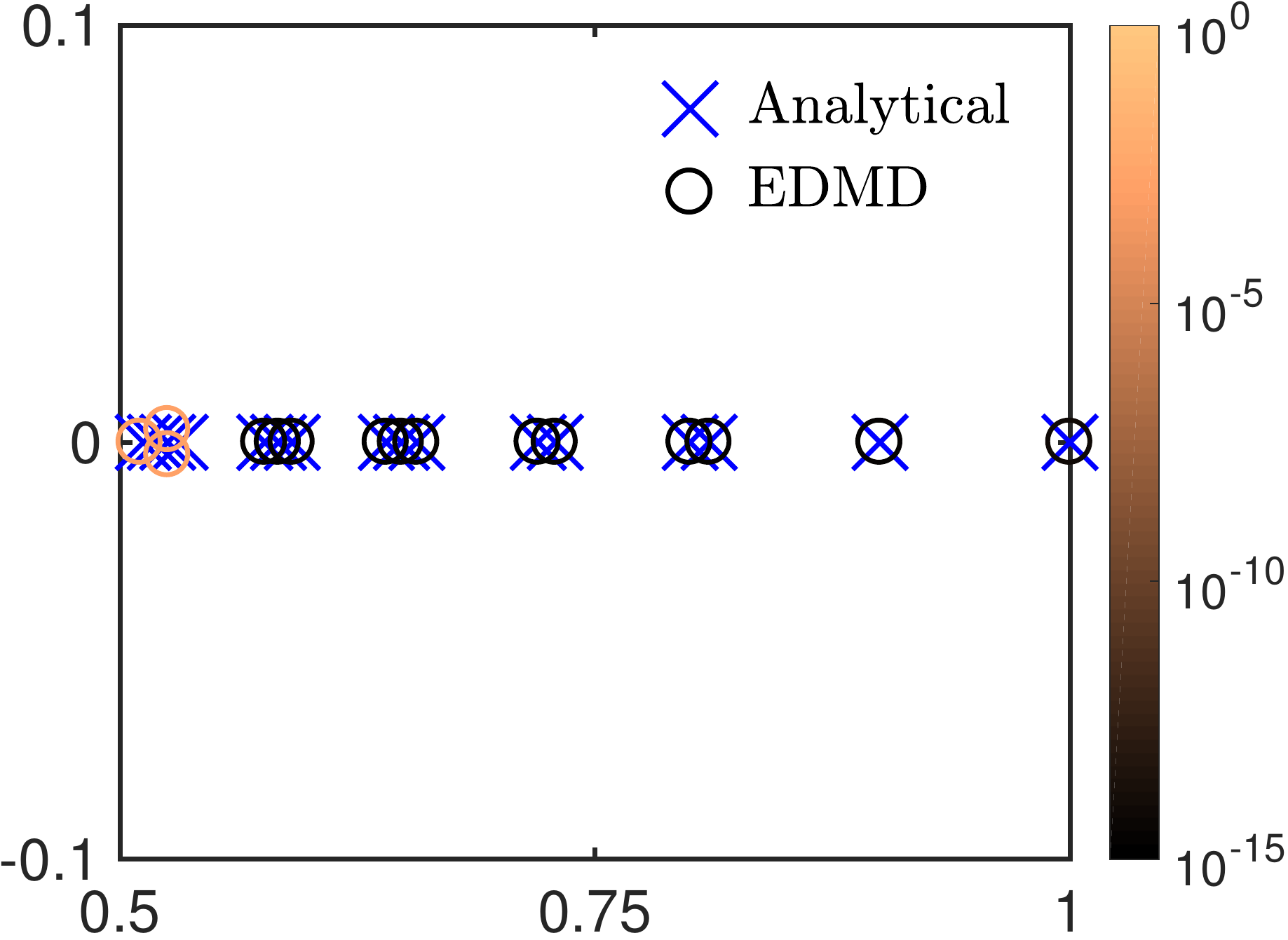}
      \small
      \node[below] at (0.5,0) {Re($\mu$)};
      \node[anchor=south,rotate=90] at (0,0.5) {Im($\mu$)};
    \end{tikzonimage}
    \caption{EDMD/analytical eigenvalues}
        \end{subfigure}%
        \begin{subfigure}[b]{0.435\textwidth}
    \centering
    \begin{tikzonimage}[width=.9\linewidth]{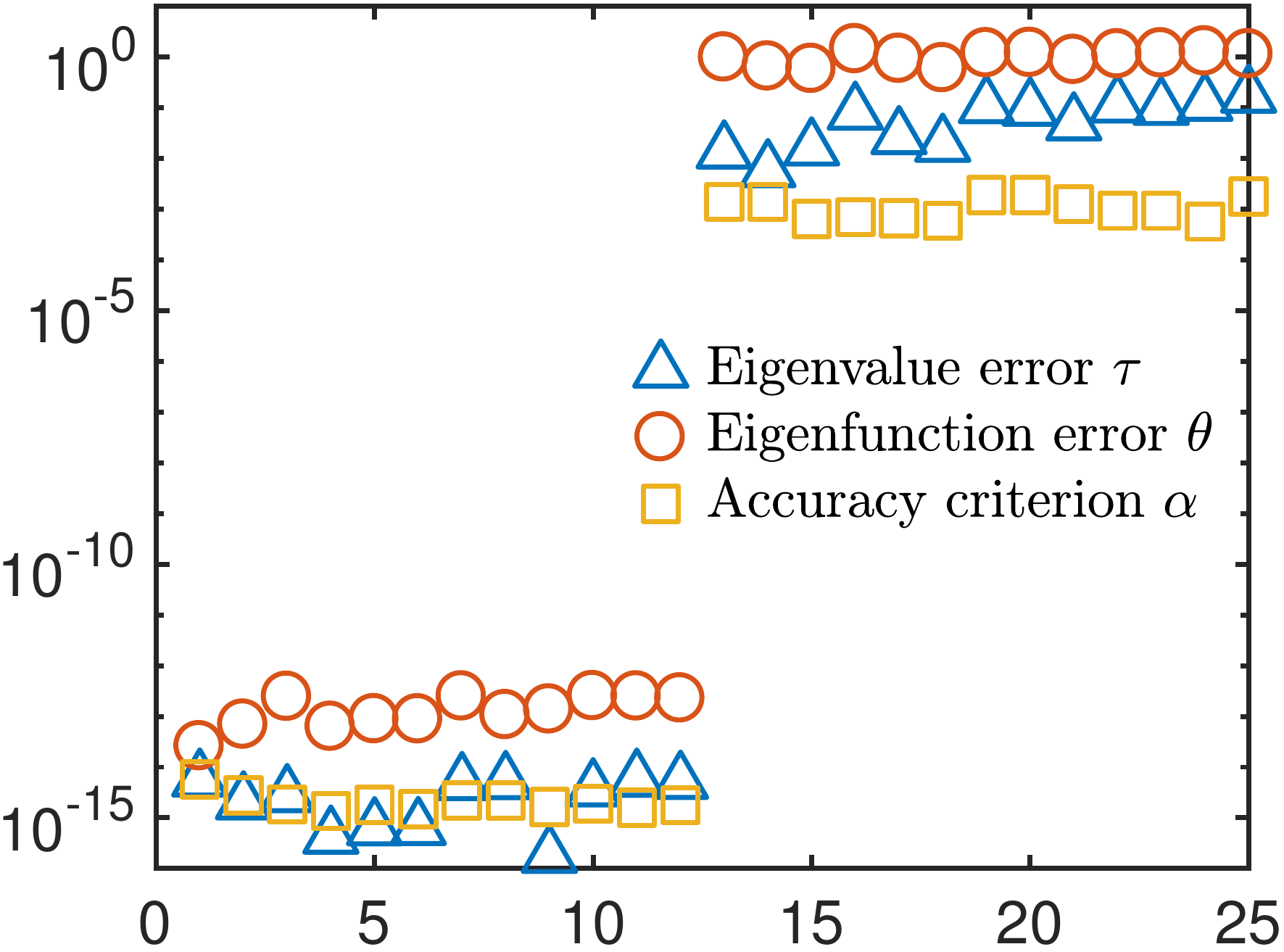}
      \small
      \node[below] at (0.5,0) {Eigenvalue index};
      \node[anchor=south,rotate=90] at (0,0.5) {Error ($\tau, \theta,\alpha$)};
    \end{tikzonimage}
    \caption{Accuracy criterion and true error}
        \end{subfigure}%
	\caption{(a) EDMD eigenvalues (circles) and analytical eigenvalues (crosses). EDMD eigenvalues are superimposed by the corresponding accuracy criterion (mode error) $\alpha$ as shown in the colorbar. (b) Comparison between the accuracy criterion $\alpha$, eigenvalue error $\tau$, and eigenfunction error $\theta$. The eigenvalues are indexed by their absolute value, in descending order.}
	\label{fig:2Dedmd}
\end{figure}

Here we apply EDMD with monomials as observables. In particular, the observables are taken to be
\begin{equation*}
\psi_{k,\ell}(\bm{x}) = x_1^k x_2^\ell, \quad k,\ell=0,1,2,3,4,5,
\end{equation*}
where the feature space dimension is $q =6 \times 6=36$.  The results of
applying EDMD are shown in Figure \ref{fig:2Dedmd}(a). We note that mode error
indicates that leading eigenvalues are approximated very accurately ($\alpha
\sim 10^{-15}$), and this is consistent with the comparison to analytical
eigenvalues. As mentioned in section~\ref{sec:DMD}, if the Koopman
eigenfunctions lie in the span of the observables, the eigenfunction can be
found exactly by EDMD. In this case, monomials up to degree~5 span the leading Koopman
eigenfunctions, and hence these eigenvalues can be identified.

To validate that the proposed accuracy criterion does indeed indicate accuracy,
now we compare $\alpha$ with the true error. We can compute the discrepancy
between DMD eigenvalues, indicated by $\hat\mu_i$, and true eigenvalues ${\mu}_{k,\ell} = \gamma^k \delta^{\ell}$ given in equation~\eqref{eq:2Danal}, by defining the eigenvalue error
\begin{equation}
\tau_i = \frac{|\hat\mu_i - {\mu}_{k,\ell}|}{|{\mu}_{k,\ell}|},
\end{equation}
where the indices $(k,\ell)$ are chosen such that ${\mu}_{k,\ell}$ is the
closest eigenvalue to $\hat \mu_i$. We then interpret $\hat \mu_i$ as a DMD approximation to the analytical eigenvalue ${\mu}_{k,\ell}$. We can also compute the discrepancy between DMD eigenfunctions $\hat\varphi_i$ and true eigenfunctions ${\varphi}_{k,\ell}$ given in equation~\eqref{eq:2Danal}. We normalize the eigenfunctions $\hat\varphi_i$ and  ${\varphi}_{k,\ell}$ so that $|\varphi|_\text{max}=1$ in the domain $\Omega = [-1,1] \times [-1,1]$, and define the eigenfunction error as
\begin{equation}
\theta_i = \frac{\|\hat\varphi_i - {\varphi}_{k,\ell}\|}{\|{\varphi}_{k,\ell}\|},
\end{equation}
where $\|\cdot\|$ denotes the $L^2$ norm given by
\begin{equation}
\|f\|^2 = \int_{\Omega} |f(\bm{x})|^2 d \bm{x}.
\end{equation}


In order to validate the accuracy criterion, we compare $\alpha_i$ with the
eigenvalue error $\tau_i$ and eigenfunction error $\theta_i$ in Figure
\ref{fig:2Dedmd}(b). We observe that $\alpha$ highly correlates with both $\tau$
and $\theta$, even though the proposed accuracy criterion does not assume access
to analytical Koopman eigenpairs. The proposed accuracy criterion hence
indicates accuracy very well, by comparison with the true error defined using true
Koopman eigenpairs. Starting from the 13th eigenvalue
$\hat{\mu}_{13}\approx \mu_{6,0}=\gamma^6 \delta^0=(0.9)^6 (0.8)^0=
0.531441$, the error dramatically increases, which implies that the remaining
eigenfunctions cannot be accurately identified using EDMD with this choice of
observables. This is expected, as monomials up to degree~5 can not represent the
eigenfunction ${\varphi}_{6,0}(\bm{x}) = x_1^6$. This comparison gives us
confidence in the reliability of the accuracy criterion.

\begin{figure}[t]
        \centering
                  \begin{tikzonimage}[width=.9\linewidth]{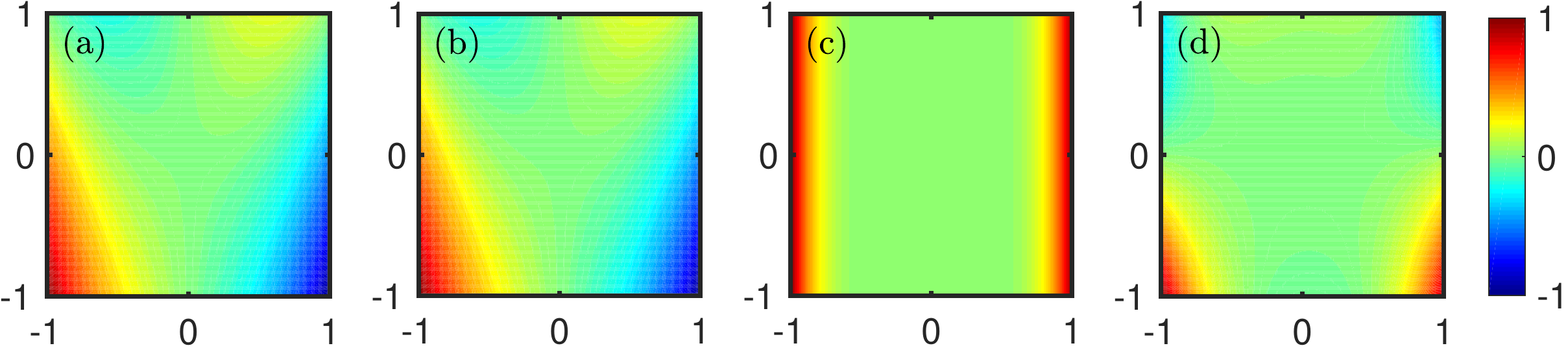}
    \small
    \node[below] at (0.125,0) {$x_1$};
    \node[above] at (0.125,1) {${\varphi}_{1,1}$};
    \node[below] at (0.36,0) {$x_1$};
        \node[above] at (0.36,1) {$\hat\varphi_6$};
    \node[below] at (0.6,0) {$x_1$};
        \node[above] at (0.6,1) {${\varphi}_{6,0}$};
    \node[below] at (0.83,0) {$x_1$};
        \node[above] at (0.83,1) {$\hat\varphi_{13}$};

    \node[anchor=south,rotate=90] at (0,0.55) {$x_2$};


  \end{tikzonimage}
  \caption{Eigenfunctions for the system defined in \eqref{eq:2Dex}, restricted
    to a domain of $[-1,1] \times [-1,1]$, and normalized such that
    $|\varphi(\bm{x})|_{max}=1$. The analytical eigenfunction ${\varphi}_{1,1}$
    shown in (a) is closely approximated by the eigenfunction
    $\hat \varphi_6$ computed by EDMD, shown in (b).  However, the
    analytical eigenfunction ${\varphi}_{6,0}$ (with eigenvalue
    $\mu_{6,0}=0.531441$) shown in (c) is not closely approximated by its
    corresponding eigenfunction $\hat \varphi_{13}$ computed by EDMD (with
    eigenvalue $\hat\mu_{13}=0.5250+0.0030j$), whose real part is shown in~(d).}
  \label{fig:EDMDefunc}
\end{figure}
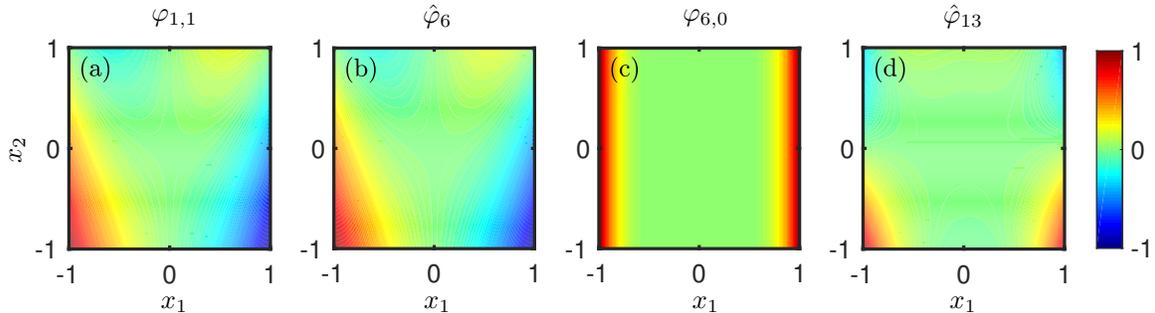

We now consider the 6th eigenvalue $\hat{\mu}_6\approx\mu_{1,1} = 0.72$,
and the 13th eigenvalue $\hat{\mu}_{13}\approx\mu_{6,0}=0.531441$. The
errors $\tau_6 \approx 10^{-15}$ and  $\theta_6 \approx 10^{-13}$ indicate that the 6th eigenpair is approximated very accurately, while $\tau_{13} \approx 10^{-2},\theta_{13} \approx 10^0$ indicate that the 13th eigenpair is approximated with lower accuracy. The EDMD eigenfunctions are compared with the analytical eigenfunctions in Figure~\ref{fig:EDMDefunc}. It is observed that the 6th eigenfunction is indeed approximated very accurately, as $\alpha_6 \approx 10^{-15}$ suggests. The 13th eigenfunction are approximated less accurately, as is expected given that $\alpha_{13} \approx 10^{-3}$. This comparison shows that the accuracy criterion does indicate the accuracy of DMD approximated Koopman eigenpairs, without assuming access to the true Koopman eigenpairs.

\section{Evaluating the performance of kernels with the accuracy criterion}
\label{sec:synth}
This section focusses on using the accuracy criterion defined in section \ref{sec:proposemetric} to evaluate the performance of KDMD using various kernel functions. We first introduce a few commonly used kernel functions in section \ref{sec:kernelfunc}, then we compare the performance of various kernels in section \ref{sec:kernelperform}, using the same test problem considered in section \ref{sec:Validating}. Following this, section \ref{sec:kernelrobust} studies  the robustness of various kernels for the case where the data are noisy.

\subsection{Kernel functions}
\label{sec:kernelfunc}
In section \ref{sec:kdmd} we briefly described  KDMD, which makes use of a kernel function to circumvent the curse of dimensionality associated with EDMD. Application of KDMD requires a suitable choice of kernel function. In order to appreciate how a kernel function may implicitly define a observable function, note that Mercer's theorem \cite{mercer1909functions} states that a (quite broad) class of ``Mercer kernels'' $k(\bm{x},\hat{\bm{x}})$ may be written as
\begin{equation}
k(\bm{x},\hat{\bm{x}})= \sum_{i=1}^{\infty} c_i \psi_i(\bm{x}) \psi_i(\hat{\bm{x}}), \quad c_i \geq c_{i+1} \geq 0.
\label{eqn:kernelexpansion}
\end{equation}
Hence there exists an infinite dimensional implicit observable function (also called feature map in the machine learning community)
\begin{equation}
\bm{\psi}(\bm{x})=
\begin{bmatrix}
\sqrt{c_1} \psi_1(\bm{x}) & \sqrt{c_2} \psi_2(\bm{x}) & \cdots & \sqrt{c_i} \psi_i(\bm{x}) & \cdots
\end{bmatrix}^T
\end{equation}
such that $k(\bm{x},\hat{\bm{x}}) = \langle \bm{\psi}(\bm{x}), \bm{\psi}(\hat{\bm{x}}) \rangle$.
We now introduce a few commonly used kernel functions, and in
section~\ref{sec:kernelperform} we compare their performance on the example from
the previous section.

\paragraph{Polynomial kernel}
\begin{equation}
k(\bm{x},\hat{\bm{x}})=(1+\bm{x}^T\hat{\bm{x}})^d
\end{equation}
The (implicit) observables associated with the polynomial kernel are all
monomials in components of $\bm{x} \in \mathbb{R}^n$ up to degree $d$. The
dimension of the observable vector is $q=\binom{n+d}{d}$. The feature map for arbitrary $n \geq 1, d \geq 0$ is described in details in \cite{cotter2011explicit}. The observables when $n=2,d=2$ are given by equation \eqref{eq:polykerneld2}.

\paragraph{Exponential kernel}
\begin{equation}
k(\bm{x},\hat{\bm{x}})=\exp\big(\bm{x}^T\hat{\bm{x}}\big)
\end{equation}
The (implicit) observables associated with the  exponential kernel are all monomials in components of $\bm{x}$, up to infinite degree. An explicit feature map can be also found from a Taylor expansion of the exponential kernel \cite{cotter2011explicit}. Taking $\bm{x} \in \mathbb{R}^2$ for example, the kernel can be expanded as
\begin{equation*}
\begin{split}
\exp\{\bm{x}^T\hat{\bm{x}}\} & =\sum_{\ell=0}^{\infty} \frac{(\bm{x}^T\hat{\bm{x}})^\ell}{\ell!}=\sum_{\ell=0}^{\infty} \frac{(x_1 \hat{x}_1+x_2 \hat{x}_2)^\ell}{\ell!} \\
& =
\sum_{\ell=0}^{\infty}  \frac{\sum_{k=0}^{\ell} \binom{l}{k} (x_1 \hat{x}_1)^k(x_2 \hat{x}_2)^{\ell-k}}{\ell!}
=\langle \bm{\psi}(\bm{x}),\bm{\psi}(\hat{\bm{x}}) \rangle,
\end{split}
\end{equation*}
where the observable is $\psi_{\ell,k}(\bm{x})=\left(\binom{\ell}{k}\big/\ell!\right)^{1/2} x_1^k x_2^{\ell-k},$ where $\ell=0,1,2,\cdots$, and $k=0,1,2,\cdots,\ell$.
Notice that the number of observables is infinite, $q=\infty$.

\paragraph{Gaussian kernel}
\begin{equation}
k(\bm{x},\hat{\bm{x}})= \exp\bigg(-\frac{\|\bm{x}-\hat{\bm{x}}\|_2^2}{\sigma^2}\bigg),
\end{equation}
where $\|\cdot\|_2$ is the $\ell^2$ norm, and $\sigma$ scales the kernel width \cite{fasshauer2011positive}.
The Gaussian kernel is a Mercer kernel for all dimensions $n \geq 1$ \cite{scovel2010radial}. Take $x \in \mathbb{R}$ as an example, the (implicit) observables as in equation \eqref{eqn:kernelexpansion} are given by \cite{gretton2013introduction}
\begin{equation*}
\psi_k(x) \propto \exp(-(d-a)x^2) H_k(x \sqrt{2d}),
\end{equation*}
where
\begin{equation*}
c_k \propto b^k, \quad b<1,
\end{equation*}
$a,b,d$ are functions of $\sigma$, and $H_k$ is the k-th order Hermite polynomial. The number of observables is infinite, $q=\infty$. For arbitrary $n$, an explicit feature map can in principle be also found from Taylor expansion of the Gaussian kernel \cite{cotter2011explicit}.

\paragraph{Laplacian kernel}
\begin{equation}
k(\bm{x},\hat{\bm{x}})=\exp\{-\frac{\|\bm{x}-\hat{\bm{x}}\|_2}{\sigma}\}
\end{equation}
Note the similarity between the Laplacian and Gaussian kernels, with the difference being that that the Laplacian kernel uses the  $\ell^2$ norm in the exponent without squaring \cite{souza2010kernel}.
For arbitrary $n$, the Laplacian kernel is a valid Mercer kernel \cite{scovel2010radial}.

\subsection{Performance of kernels}
\label{sec:kernelperform}
\begin{figure}[t]
        \centering
                  \begin{tikzonimage}[width=.9\linewidth]{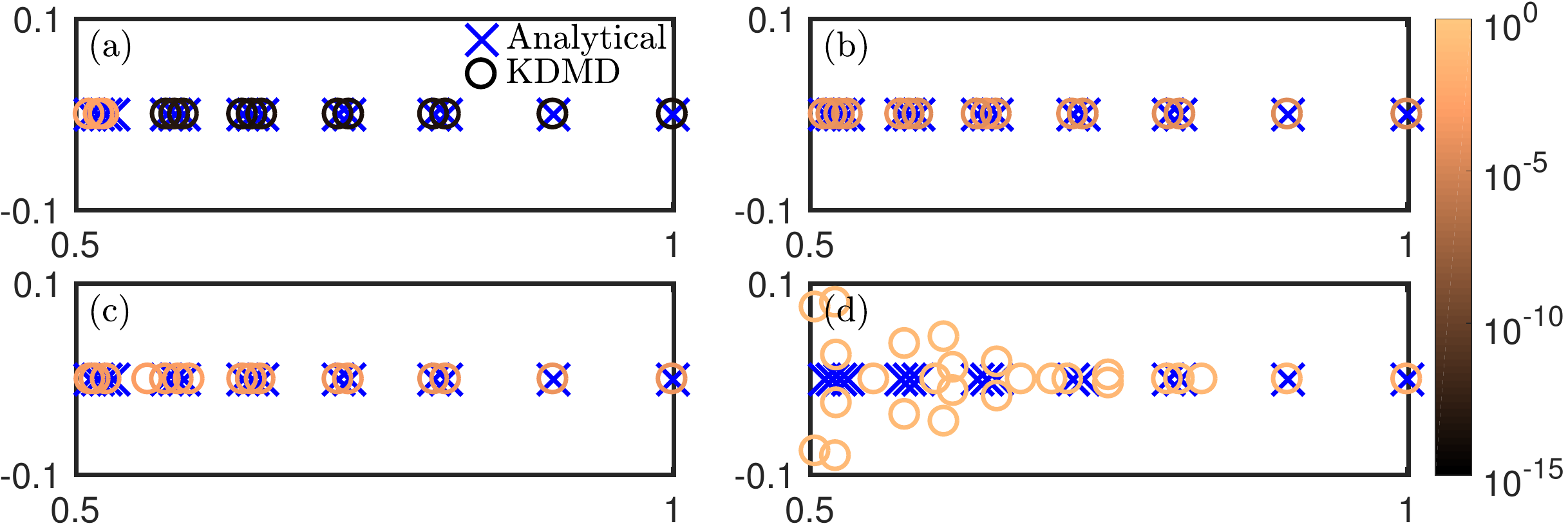}
    \small
    \node[below] at (0.23,0) {Re($\mu$)};
    \node[below] at (0.7,0) {Re($\mu$)};
    \node[anchor=south,rotate=90] at (0,0.25) {Im($\mu$)};
    \node[anchor=south,rotate=90] at (0,0.75) {Im($\mu$)};
  \end{tikzonimage}

  \caption{KDMD eigenvalues (circles) colored by their estimated mode error $\alpha$. Analytical eigenvalues (crosses) are shown for comparison. (a) Polynomial kernel of degree $d=5$, $q=\binom{2+5}{5}=21$. (b) Exponential kernel,  $q=\infty$. (c) Gaussian kernel with $\sigma=1$,  $q=\infty$. (d) Laplacian kernel with $\sigma=1$,  $q=\infty$.}
  \label{fig:KDMD1}
\end{figure}
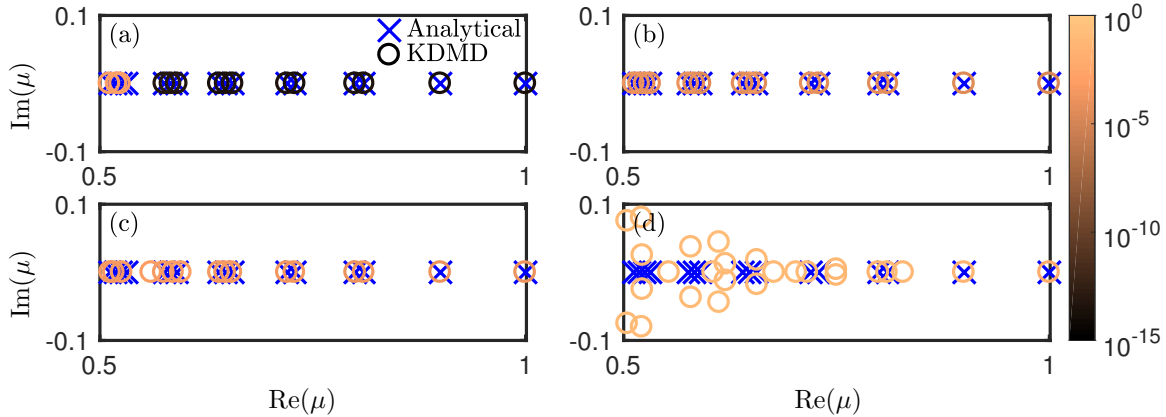

We now compare the above kernel functions using the example considered in
section~\ref{sec:Validating}.  Figure~\ref{fig:KDMD1} shows the performance of
polynomial, exponential, Gaussian, and Laplacian kernels in identifying the
Koopman eigenvalues of the system, using the same training and testing data as
in section~\ref{sec:Validating}.

We find that a polynomial kernel of degree $d=5$ accurately identifies the
leading eigenvalues (${\mu}_{k,\ell} \in [0.6,1]$) with very high accuracy
($\alpha \approx 10^{-14}$), as was the case with EDMD. This is not surprising,
as the polynomial kernel implicitly defines monomials of states as observables,
which span the same space as the explicitly defined monomials used in EDMD. With
increasing order of the polynomial kernel, more eigenvalues can be accurately
identified. It is found that the exponential kernel can identify more
eigenvalues (${\mu}_{k,\ell} \in [0.5,1]$) than the polynomial kernel with
satisfactory accuracy ($\alpha \approx 10^{-4}$), since the implicit observables
associated with exponential kernel are monomials up to infinite degree. The
Gaussian kernel is able to find the leading eigenvalues (${\mu}_{k,\ell} \in
[0.65,1]$) with mode error (accuracy criterion) $\alpha \approx 10^{-4}$
to~$10^{-3}$, even though the implicit observables of the Gaussian kernel are
not monomials. This demonstrates the potential power of kernel functions: they are able to
span a useful function space, primarily because the dimension of the space of
(implicit) observables can be large, and even infinite. The Laplacian kernel can approximate only a few leading eigenvalues ($\mu=1.0.9,0.8$), and with a lower accuracy of $\alpha \approx 10^{-2}$.

We  emphasize that, while the exact Koopman eigenvalues are known in this case, it is possible to use the accuracy criterion to compare the performance of different kernels even when the true dynamics are unknown. Indeed, using only the results of the accuracy criterion, we would reason that the polynomial kernel is the best choice for identifying the leading Koopman eigenvalues accurately.

\subsection{Sensitivity of kernels to noise}
\label{sec:kernelrobust}
In practice, data is typically corrupted with noise. Here we present a study of
the sensitivity of different kernels with respect to the presence of noise. We
add zero-mean Gaussian noise with standard deviation
$\sigma_\text{noise}=10^{-3}$ to the 100 random uniformly distributed data pairs
taken from $[-1,1] \times [-1,1]$. The training data is noisy, but the testing
data is ``clean''. Therefore, the accuracy criterion only accounts for the
accuracy of DMD approximated Koopman eigenpairs.

\begin{figure}[t]
        \centering
                  \begin{tikzonimage}[width=.9\linewidth]{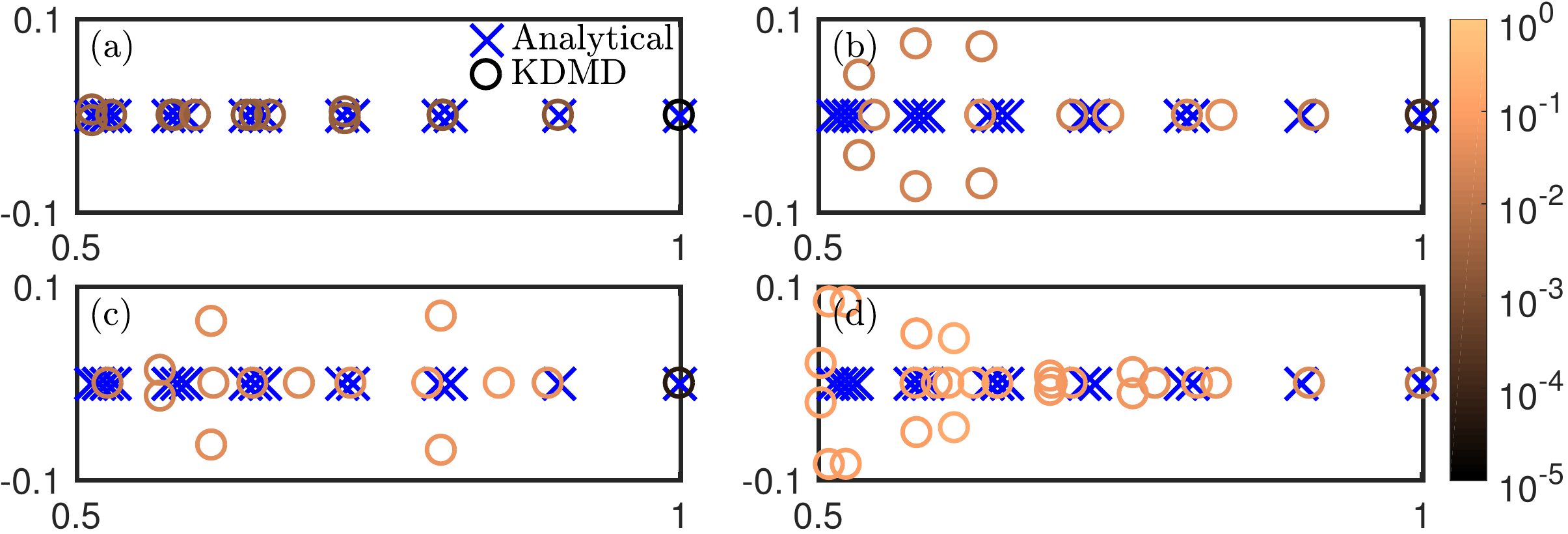}
    \small
    \node[below] at (0.23,0) {Re($\mu$)};
    \node[below] at (0.7,0) {Re($\mu$)};
    \node[anchor=south,rotate=90] at (0,0.25) {Im($\mu$)};
    \node[anchor=south,rotate=90] at (0,0.75) {Im($\mu$)};
  \end{tikzonimage}
	\caption{KDMD eigenvalues (circles) colored by their estimated mode error $\alpha$, identified from noisy data. Analytical eigenvalues (crosses) are shown for comparison. (a) Polynomial kernel of degree $d=5$, $q=\binom{2+5}{5}=21$. (b) Exponential kernel,  $q=\infty$. (c) Gaussian kernel with $\sigma=1$,  $q=\infty$. (d) Laplacian kernel with $\sigma=1$,  $q=\infty$.}
	\label{fig:KDMD2}
\end{figure}
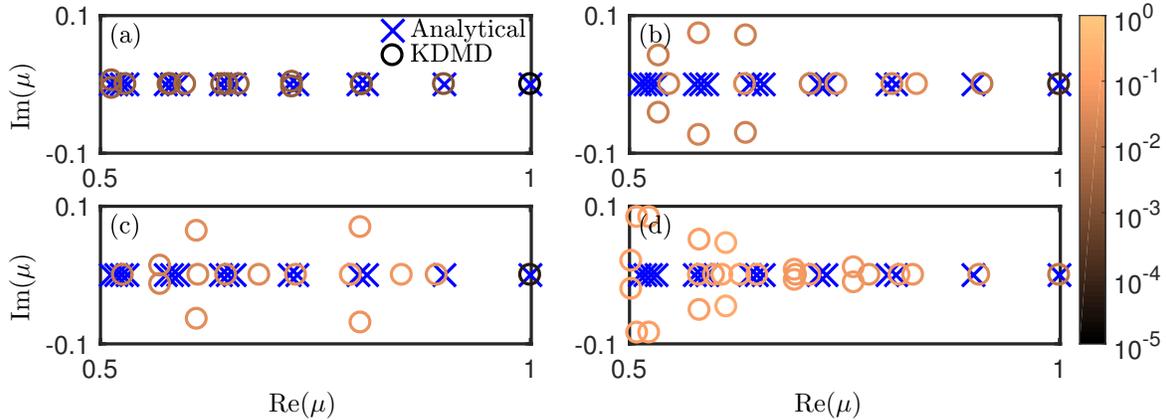

The results are shown in Figure~\ref{fig:KDMD2}. We observe that the polynomial
kernel is slightly more robust  than the other kernels ($\alpha \approx
10^{-3}$) in the presence of noise, and is able to accurately identify the first
few leading eigenvalues ($\mu=1,0.9$). The reason for this is that the dimension
of the implicit observables associated with the polynomial kernel is finite and
small ($q=21$) in comparison to the number of snapshots ($m=100$), so we avoid
problems of overfitting. In KDMD, the
Koopman eigenpairs are found from the eigendecomposition of the matrix
$\bm{A}_{\text{KDMD}}=\bm{Y}^+\bm{Y}^\#$, where the columns of $\bm{Y}$ and $\bm{Y}^\#$
are $\bm{y}=\bm{\psi}(\bm{x}) \in \mathbb{R}^q$ and
$\bm{y}^{\#}=\bm{\psi}(\bm{x}^{\#}) \in \mathbb{R}^q$ respectively, and $\bm{Y},
\bm{Y}^\# \in \mathbb{R}^{q \times m}$. The matrix~$\bm{A}_{\text{KDMD}}$ has the same non-zero
eigenvalues as the DMD matrix $\bm{A}=\bm{Y}^\# \bm{Y}^+$. $\bm{A}$ is the
optimal (least-square or minimum-norm) solution to $\min_{\bm{A}}
\|\bm{A}\bm{Y}-\bm{Y}^\#\|_F$, where $\bm{Y}, \bm{Y}^\# \in \mathbb{R}^{q \times
  m}$. For the polynomial kernel, $\bm{A}$ is the solution to an
over-constrained problem ($q<m$), and is hence more robust to noise. In
contrast, the exponential kernel, Gaussian kernel, and Laplacian kernel span an
infinite dimensional space of observables ($q=\infty$). The finite dimensional
approximation to the Koopman operator is found by solving an under-constrained
problem ($q\gg m$), which makes it more sensitive to noise, as these three
kernels tend to over-fit the noise in the trainning dataset. Given noisy data,
they are only able to accurately identify the eigenvalue $\mu=1$, whose
eigenfunction is a constant.

\section{Identifying accurate DMD modes using experimental data}
\label{sec:exp}
Having demonstrated the use of the accuracy criterion with synthetic data, now
we turn our attention to  data from fluids  experiments. In these cases, the
analytical Koopman spectral decomposition is unknown. An important advantage of
the proposed accuracy criterion is that it does not rely on known Koopman
eigenpairs, and can be applied so long as there is data available. We will use
the proposed accuracy criterion to identify accurate DMD modes for vorticity
data from flow past a circular cylinder in section \ref{sec:cylinder}, and from
a separation experiment in section \ref{sec:separation}.

\subsection{Flow past a circular cylinder}
\label{sec:cylinder}
\begin{figure}[thb]
        \centering
        \begin{subfigure}[b]{0.45\textwidth}
            \centering
    \begin{tikzonimage}[width=0.9\linewidth]{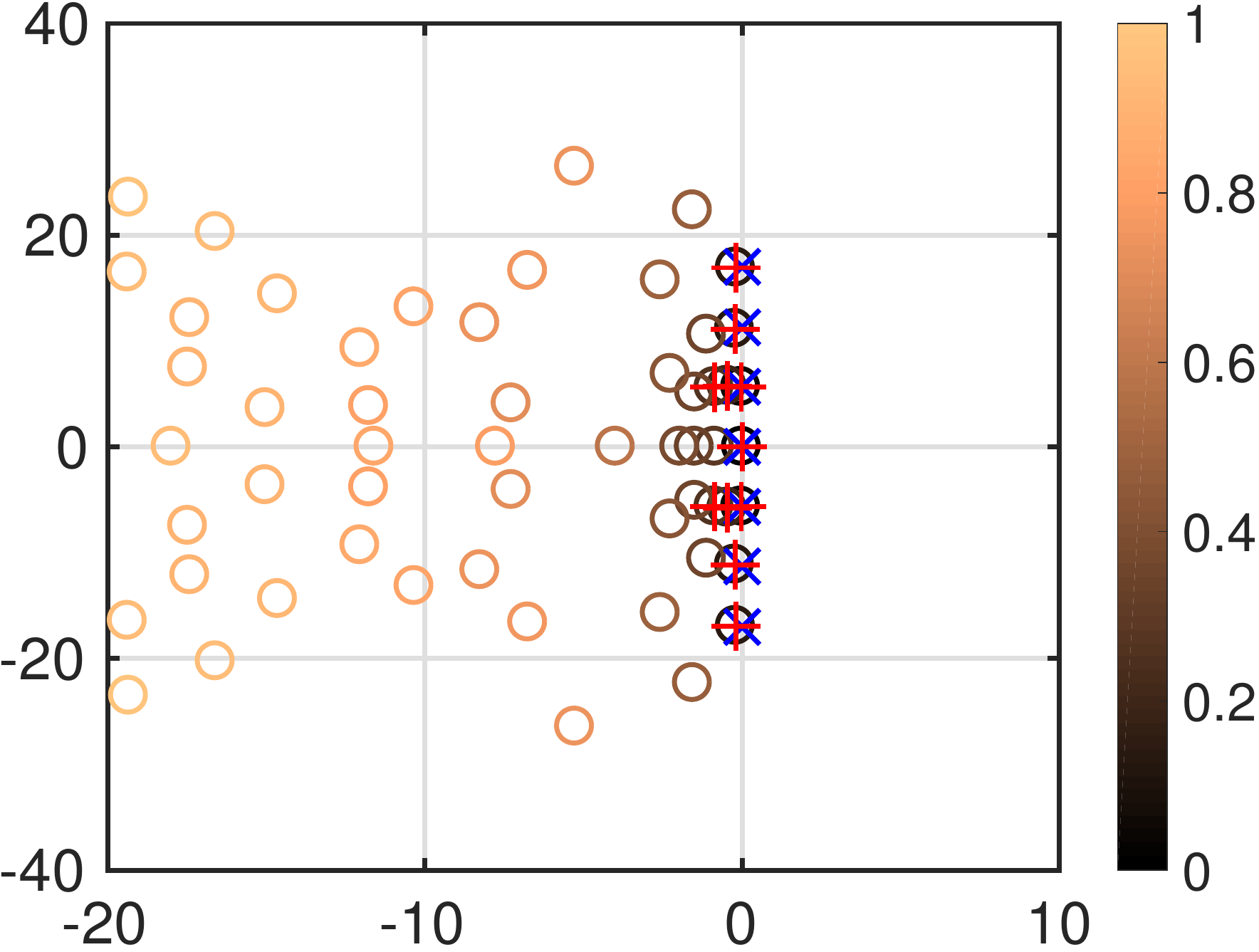}
      \small
      \node[below] at (0.5,0) {Re($\lambda$)};
      \node[anchor=south,rotate=90] at (0,0.5) {Im($\lambda$)};
    \end{tikzonimage}
        \caption{DMD, accuracy criterion}
        \end{subfigure}%
        \begin{subfigure}[b]{0.46\textwidth}
            \centering
    \begin{tikzonimage}[width=0.9\linewidth]{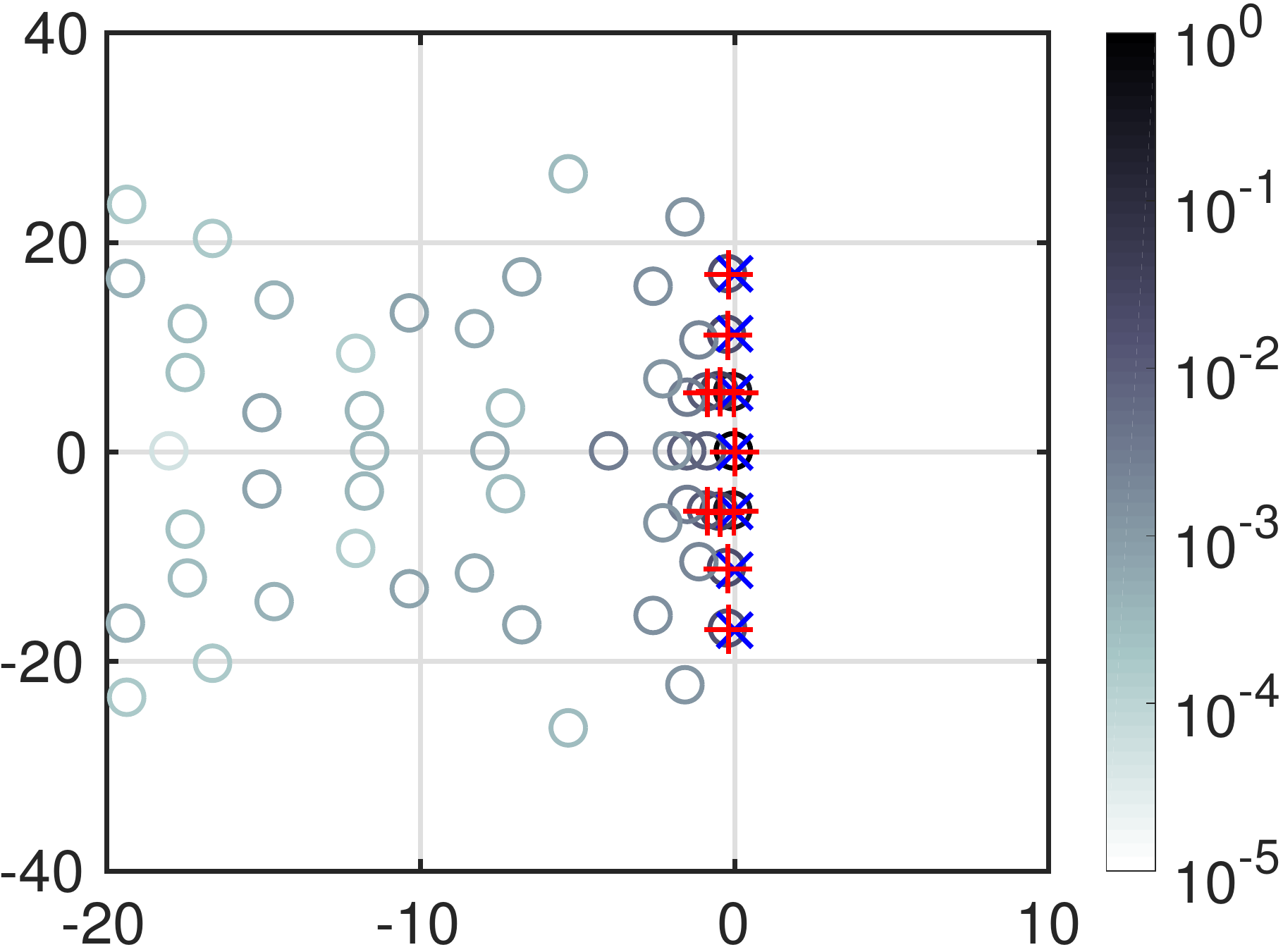}
      \small
      \node[below] at (0.5,0) {Re($\lambda$)};
      \node[anchor=south,rotate=90] at (0,0.5) {Im($\lambda$)};
    \end{tikzonimage}
                \caption{DMD, mode amplitude}
        \end{subfigure}%
        \\
          \begin{subfigure}[b]{0.3\textwidth}
    \includegraphics[width=0.95\linewidth]{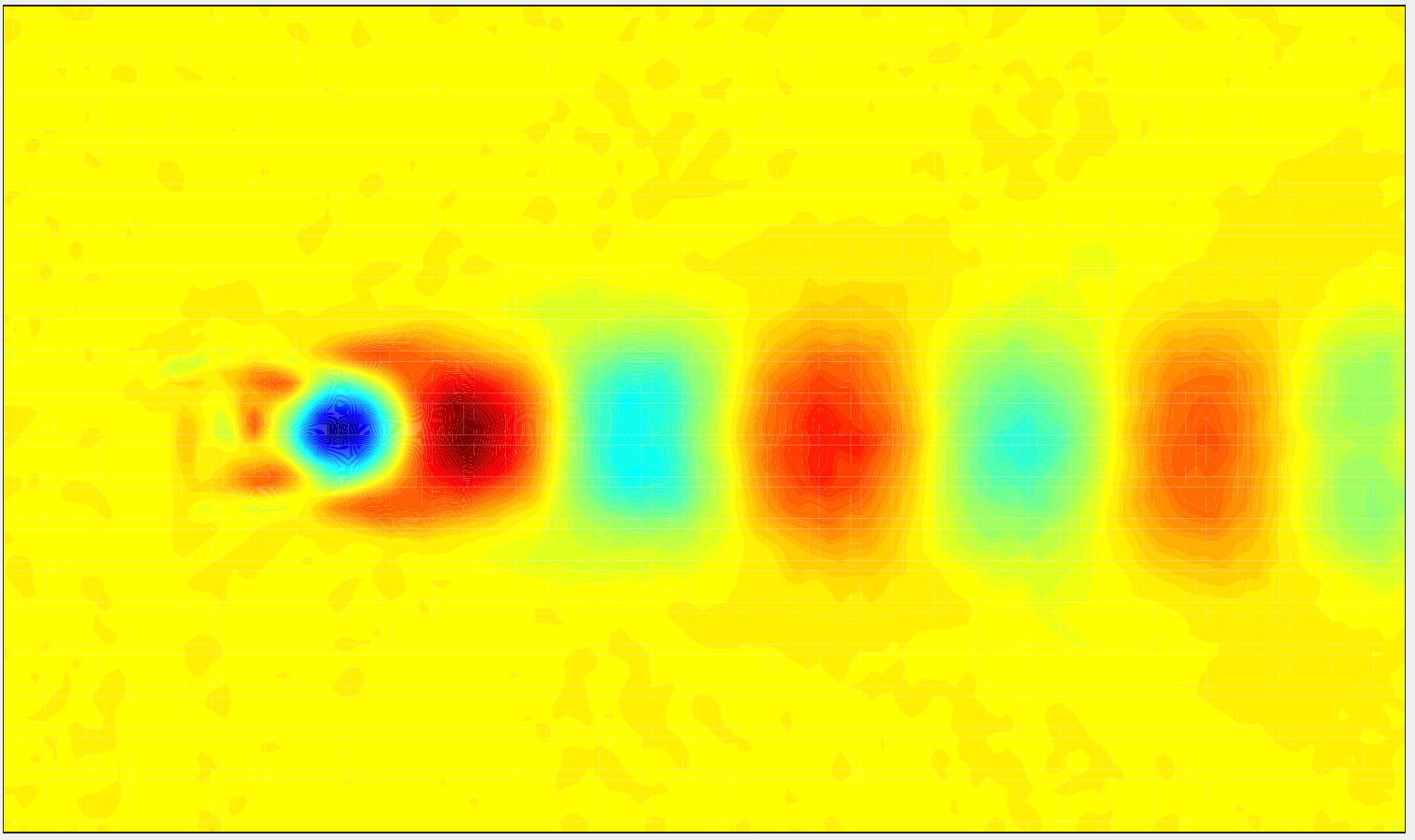}
        \centering
    \caption{Mode 1, $f_1=0.90$ Hz}
  \end{subfigure}%
  \begin{subfigure}[b]{0.3\textwidth}
    \includegraphics[width=0.95\linewidth]{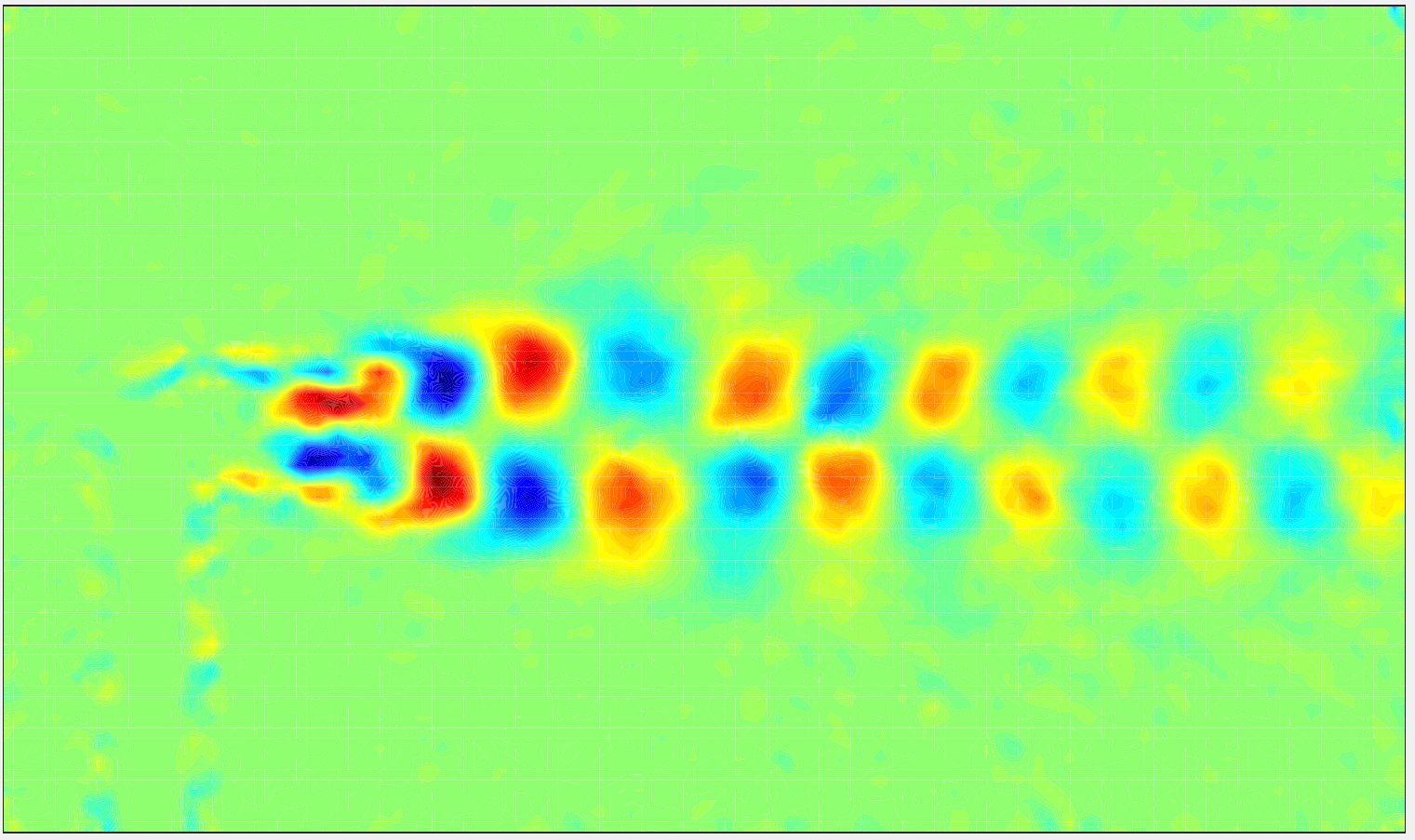}
        \centering
    \caption{Mode 2, $f_2=1.77$ Hz}
  \end{subfigure}%
    \begin{subfigure}[b]{0.3\textwidth}
    \includegraphics[width=0.95\linewidth]{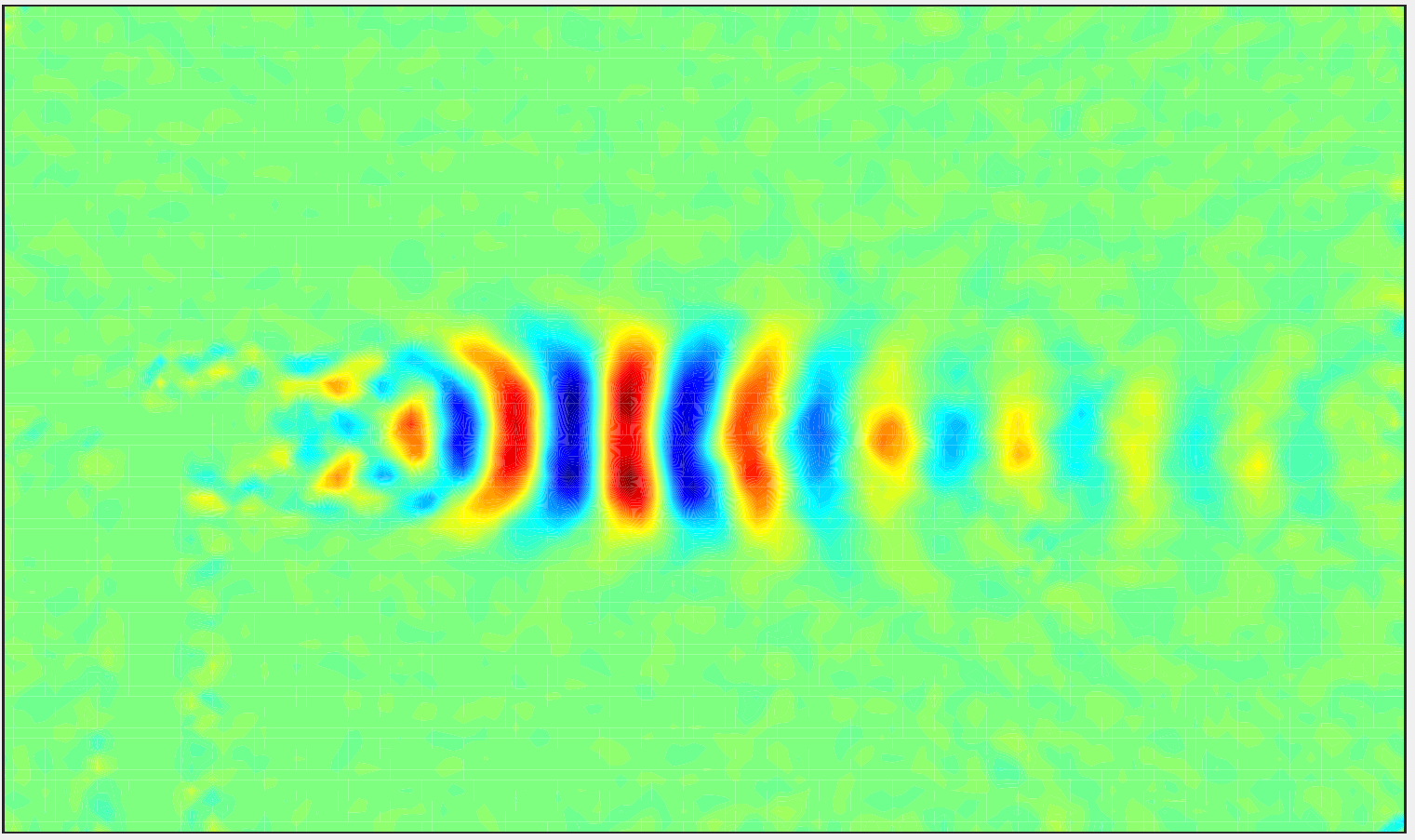}
        \centering
    \caption{Mode 3, $f_3=2.69$ Hz}
  \end{subfigure}%

	\caption{(a)-(b), Continuous-time DMD eigenvalues (circles) colored by (a) the accuracy criterion $\alpha$ and (b) mode amplitude $\beta$. Mode amplitudes are normalized by the maximum amplitude. Dominant frequencies (blue cross sign $\times$) are shown for comparison. The first 11 eigenvalues that have small $\alpha$ and large $\beta$ are shown (red plus sign $+$). (c)-(e) Three dominant DMD modes (only show real part) picked out by accuracy criterion and mode amplitude. }
	\label{fig:Re413DMD}
\end{figure}

In this example, we use the experimental particle image velocimetry (PIV) data for flow past a circular cylinder at a Reynolds number of~413. The PIV velocity data was sampled at frequency of 20\,Hz with a resolution of $135 \times 80$ pixels. See \cite{tu2014spectral} for more details about this experiment. This dataset has been used in other studies \cite{hemati2014dynamic,williams2014kernel} for testing various proposed DMD algorithms. We will use  vorticity data for DMD, which can be computed from velocity data by finite difference methods. The state dimension is $n=135 \times 80=10800$, and the number of snapshots in training data is taken to be $m = 1000$. We use an additional $m_{\text{test}}=1000$ snapshot pairs as testing data.

When we apply DMD to sequential data that has time step $\triangle t$, the continuous-time DMD eigenvalues $\lambda_{\text{DMD}}$ are related to the discrete-time DMD eigenvalues $\mu_{\text{DMD}}$ by
\begin{equation}
\mu_{\text{DMD}} = e^{\lambda_{\text{DMD}} \triangle t}.
\label{eq:cevals}
\end{equation}
The discrete-time DMD eigenvalues are computed with DMD and converted to continuous-time DMD eigenvalues by equation \eqref{eq:cevals}, and in this example the time spacing is $\triangle t=(1/20)s$. The DMD frequency $f_{\text{DMD}}$ is related to the continuous-time DMD eigenvalues $\lambda_{\text{DMD}}$ by
\begin{equation}
f_{\text{DMD}} = \frac{\text{Im}(\lambda_{\text{DMD}})}{2 \pi},
\end{equation}
where $\text{Im}(\lambda_{\text{DMD}})$ is the imaginary part of $\lambda_{\text{DMD}}$.

We first apply the standard DMD method described in section~\ref{sec:DMD}. We use a truncation level $r=100$, which corresponds to preserving $78.16\%$ of the total energy of the snapshots. 
 The continuous-time DMD eigenvalues are shown shaded by the corresponding accuracy
 criterion values~$\alpha$ in Figure~\ref{fig:Re413DMD}(a), and time-averaged
 mode amplitudes~$\beta$ in Figure~\ref{fig:Re413DMD}(b) (defined as in~\cite{kou2016improved}).

Inspecting Figure \ref{fig:Re413DMD} (a), we observe that eigenvalues near the imaginary axis are more accurate, and this observation is consistent with physical intuition: this flow exhibits a von K\'{a}rm\'{a}n vortex street, whose dominant dynamics evolve on a limit cycle. For this experiment, the wake shedding frequency is $f_\text{wake}=0.889$\,Hz \cite{tu2014spectral},
In previous work \cite{tu2014spectral}, the physically relevant dominant frequencies are reported as $f_0=0$\,Hz, $f_1=0.89$\,Hz, $f_2=1.77$\,Hz, $f_3=2.73$\,Hz.
The DMD mode associated with $\lambda_0$ is the mean of the flow, and $\lambda_1, \ \lambda_2$ and $\lambda_3$ are the first, second, and third harmonic of the fundamental wake frequency $\lambda_\text{wake}$. These four frequencies represent the dominant dynamics in this flow. This observation indicates that the proposed accuracy criterion can be used to identify physically relevant DMD modes/eigenvalues, and distinguish relevant modes from irrelevant ones. By comparing Figure~\ref{fig:Re413DMD}(a) and (b), we verify that the accuracy criterion indicates the same dominant frequencies as the mode amplitude. The DMD modes that have higher accuracy as indicated by the accuracy criterion are shown in Figure~\ref{fig:Re413DMD}(c)-(e). We verify that they look similar to those identified in previous work \cite{tu2014spectral}.



\paragraph{KDMD}
\begin{figure}[thb]
  \centering
\begin{subfigure}[b]{0.45\textwidth}
            \centering
    \begin{tikzonimage}[width=0.9\linewidth]{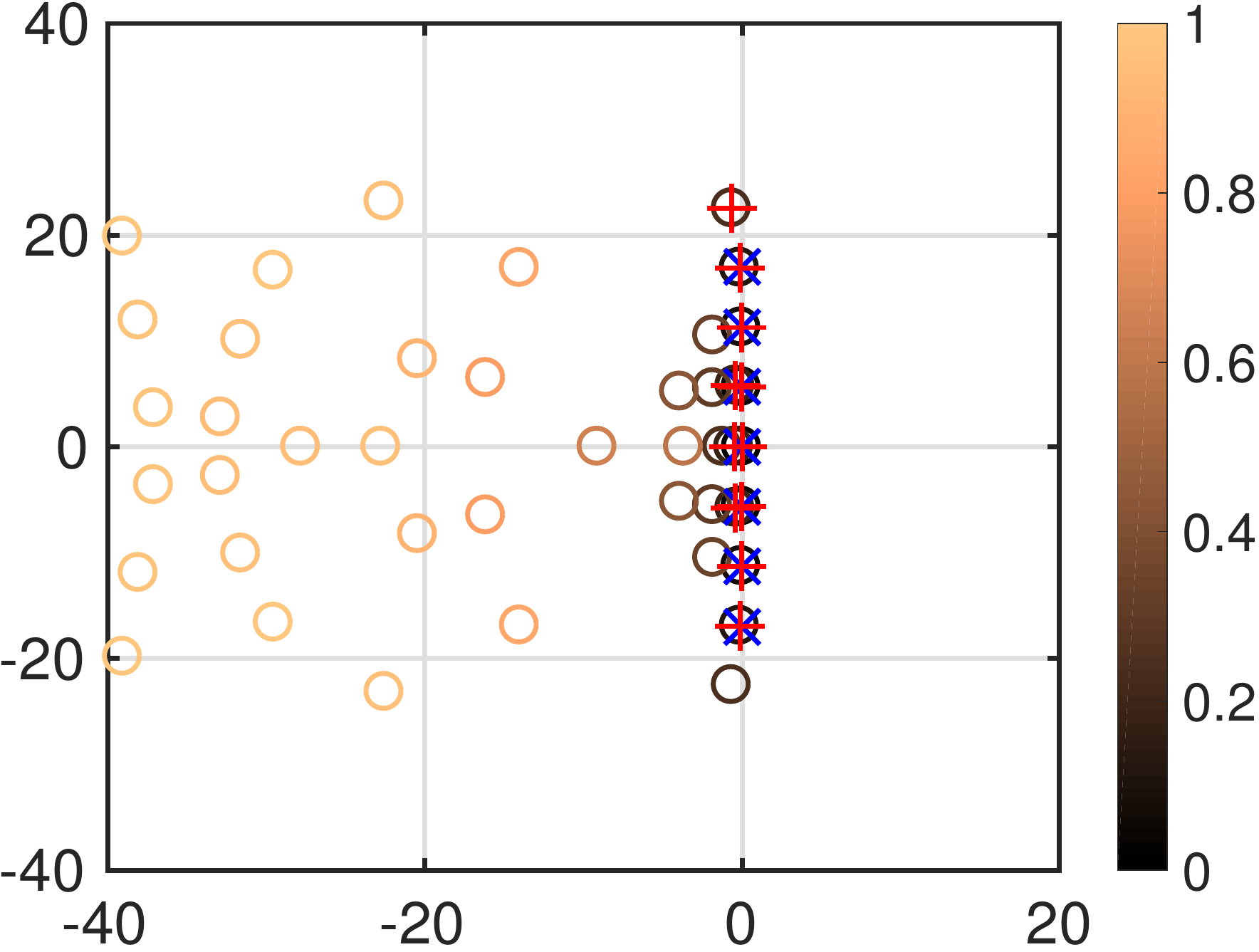}
      \small
      \node[below] at (0.5,0) {Re($\lambda$)};
      \node[anchor=south,rotate=90] at (0,0.5) {Im($\lambda$)};
    \end{tikzonimage}
        \caption{KDMD, accuracy criterion}
        \end{subfigure}%
        \begin{subfigure}[b]{0.46\textwidth}
            \centering
    \begin{tikzonimage}[width=0.9\linewidth]{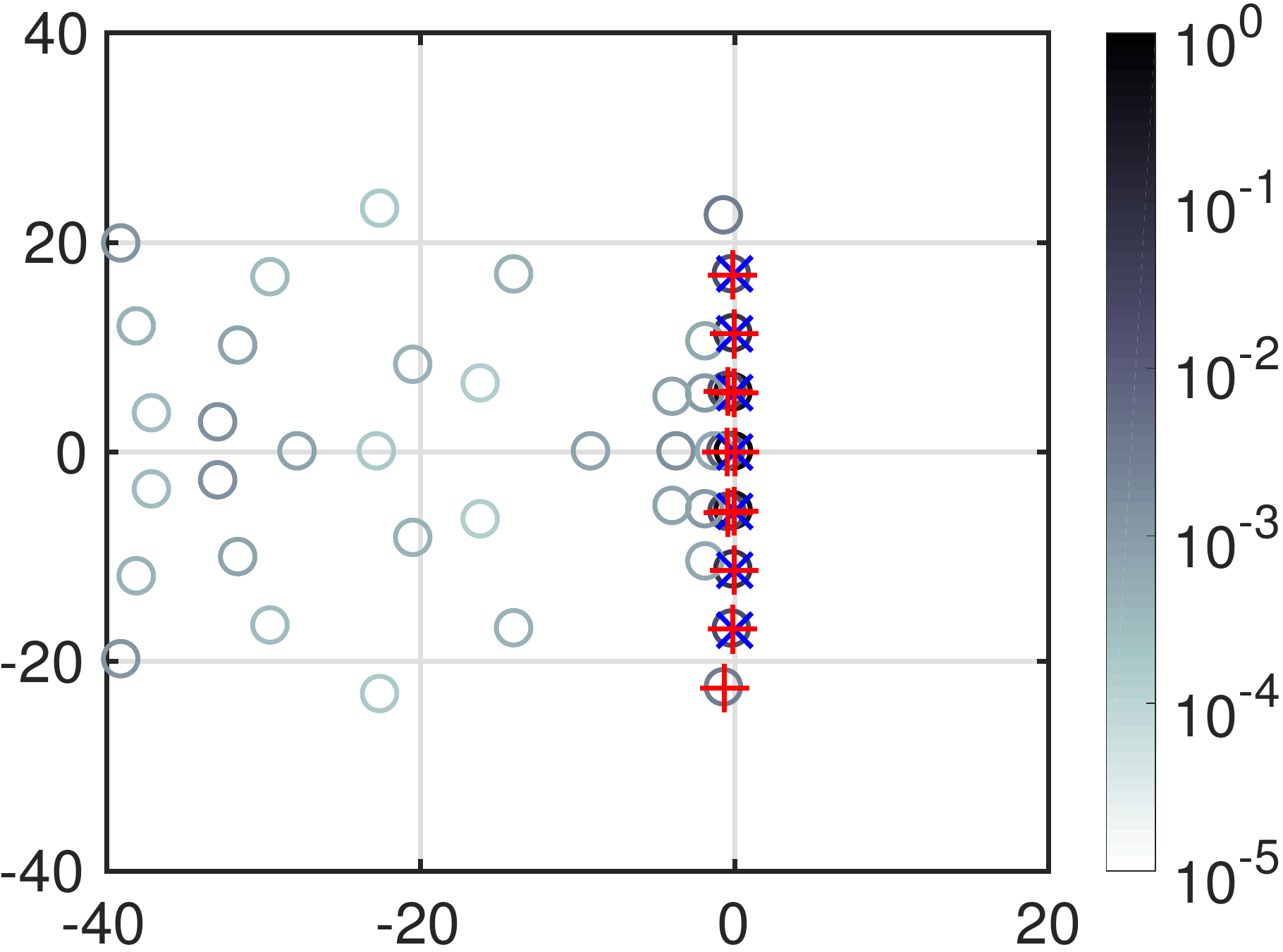}
      \small
      \node[below] at (0.5,0) {Re($\lambda$)};
      \node[anchor=south,rotate=90] at (0,0.5) {Im($\lambda$)};
    \end{tikzonimage}
                \caption{KDMD, mode amplitude}
        \end{subfigure}%
  \\
    \begin{subfigure}[b]{0.3\textwidth}
    \includegraphics[width=0.95\linewidth]{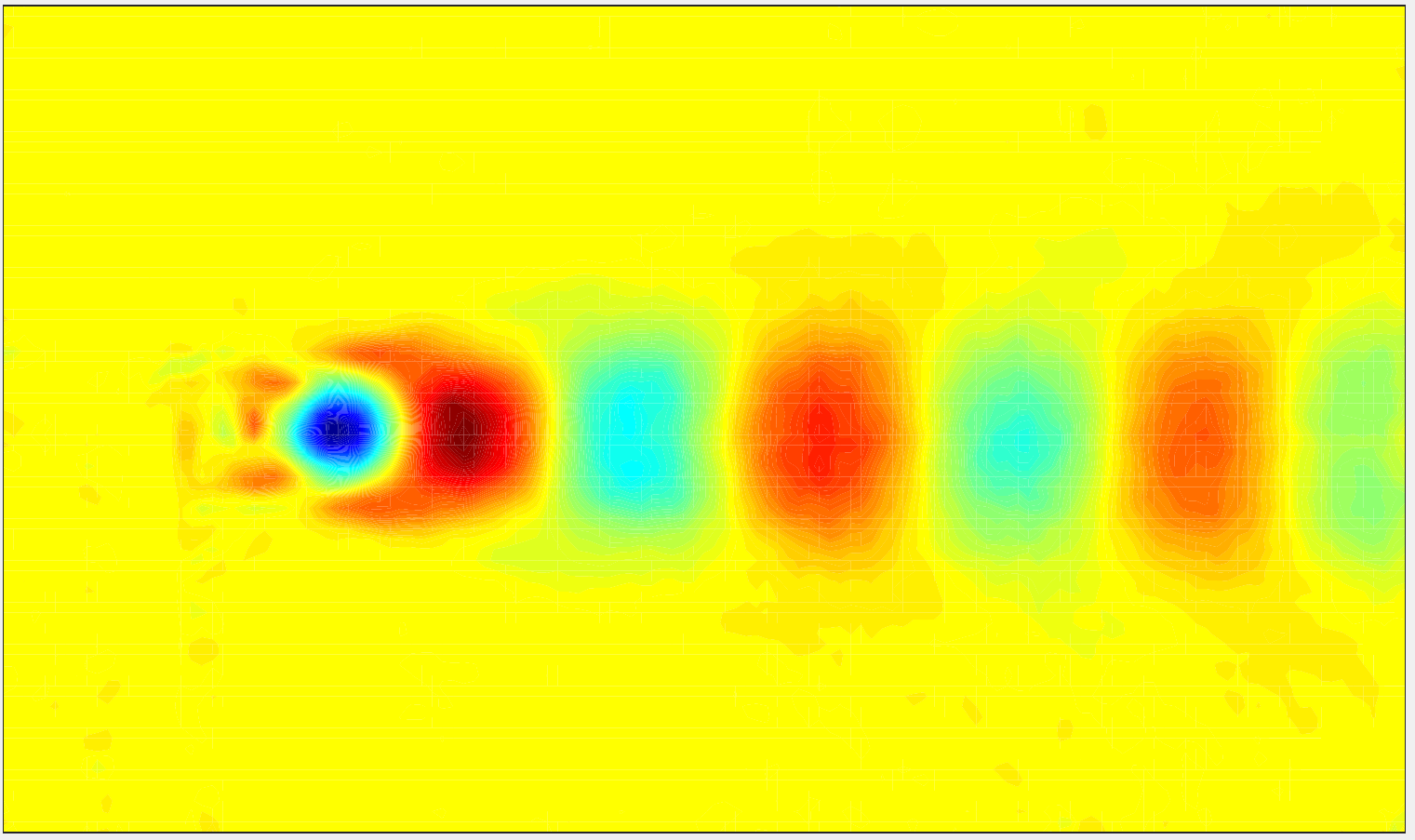}
        \centering
    \caption{Mode 1, $f_1=0.90$ Hz}
  \end{subfigure}%
  \begin{subfigure}[b]{0.3\textwidth}
    \includegraphics[width=0.95\linewidth]{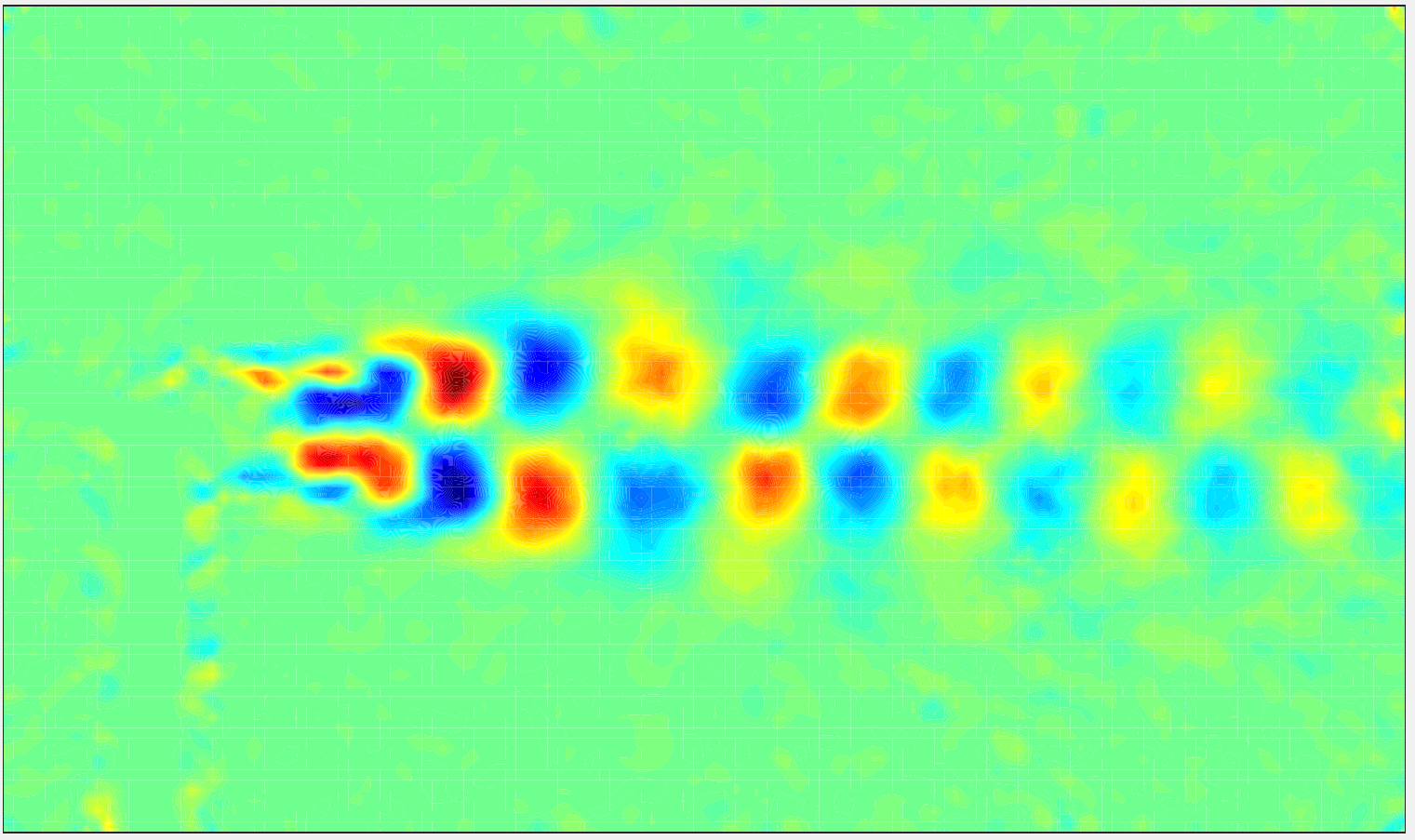}
        \centering
    \caption{Mode 2, $f_2=1.79$ Hz}
  \end{subfigure}%
    \begin{subfigure}[b]{0.3\textwidth}
    \includegraphics[width=0.95\linewidth]{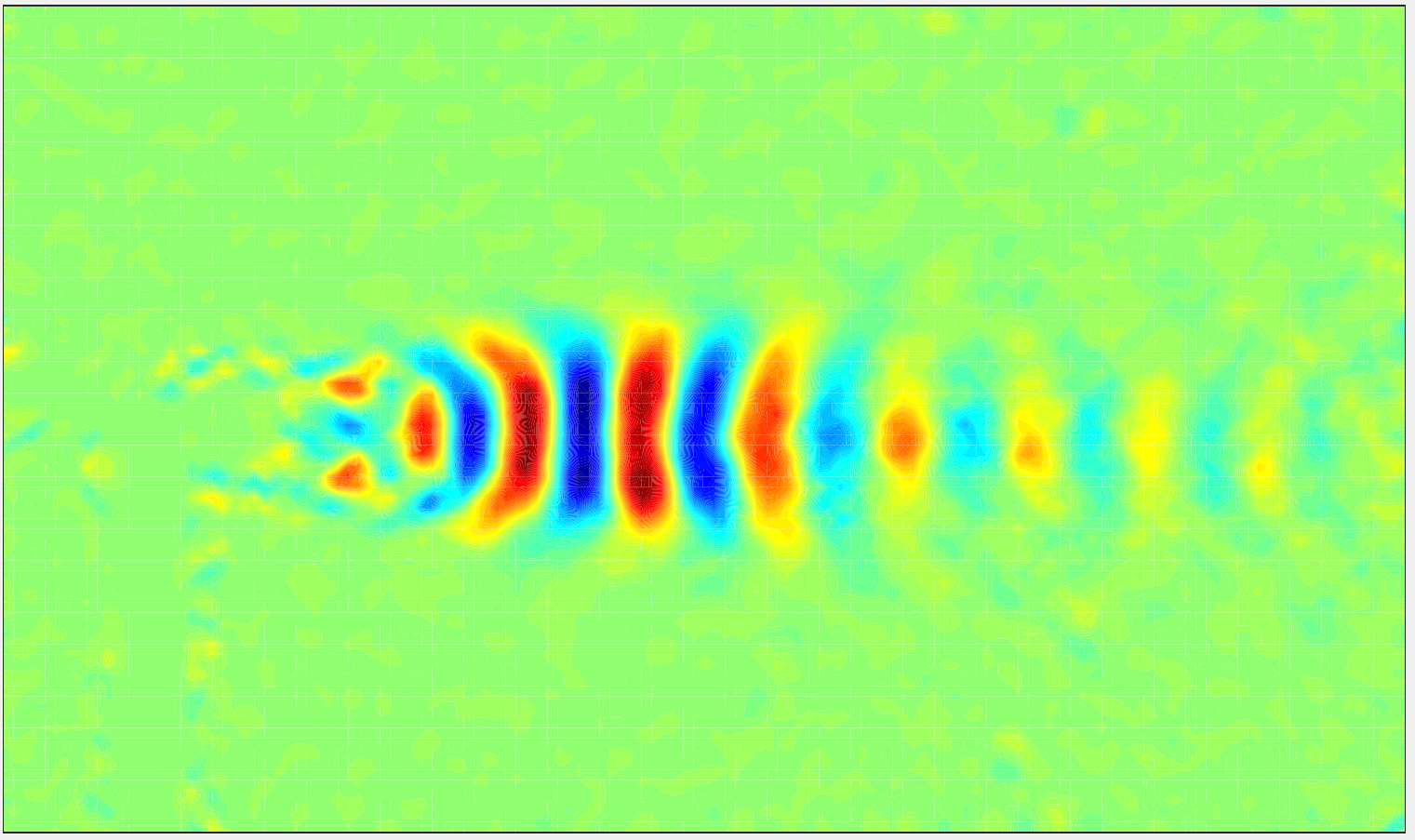}
        \centering
    \caption{Mode 3, $f_3=2.69$ Hz}
  \end{subfigure}%
  \caption{(a)--(b), Continuous-time KDMD eigenvalues (circles) colored by (a) the accuracy criterion $\alpha$ and (b) mode amplitude $\beta$.  Mode amplitudes are normalized by the maximum amplitude. Dominant frequencies (blue cross sign $\times$) are shown for comparison. The first 11 eigenvalues that have small $\alpha$ and large $\beta$ are shown (red plus sign $+$). (c)--(e) Three dominant DMD modes (real part) picked out by accuracy criterion and mode amplitude. }
  \label{fig:Re413KDMD}
\end{figure}

Next, we investigate the performance of KDMD on this dataset.
Figure~\ref{fig:Re413KDMD} shows results for a polynomial kernel of degree $d=5$, again using a truncation level of $r=100$. The DMD eigenvalues are shown in Figure~\ref{fig:Re413KDMD}(a)--(b), colored by both accuracy criterion and mode amplitude. The relevant DMD modes picked out by accuracy criterion and mode amplitude are shown in Figure~\ref{fig:Re413KDMD}(c)--(e). We verify that accuracy criterion is able to isolate dominant modes when using KDMD.

\subsection{Canonical separated flow}
\label{sec:separation}
\begin{figure}[thb]
        \centering
        \begin{subfigure}[b]{0.5\textwidth}
        \includegraphics[width=0.95\linewidth]{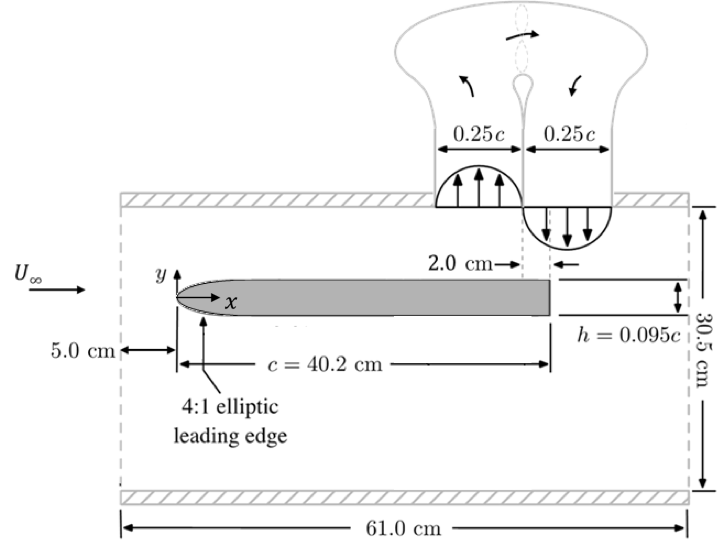}
                \centering
        \caption{Experiment setup}
        \end{subfigure}%
        \begin{subfigure}[b]{0.5\textwidth}
                \includegraphics[width=0.95\linewidth]{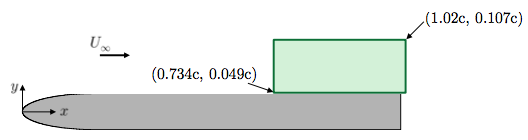}
                        \centering
                \caption{The PIV measurement region}
        \end{subfigure}%
	\caption{Sketch of the canonical separated flow experiment setup (adapted from \cite{griffin2013control}) and the PIV measurement region.}
	\label{fig:separation}
\end{figure}

In this example, we use PIV data from a canonical flow separation experiment
sketched in Figure~\ref{fig:separation}. Separation is induced on the surface of
a flat plate by a suction/blowing boundary condition imposed on the wall of the
wind tunnel, near the trailing edge of the plate. The free-stream velocity is
$U_{\infty}=3.9\,\text{m/s}$, the chord length is $c=402$\,mm, the span is $s=305$\,mm, and
the height is $h=0.095c$. The Reynolds number based on chord length is
$Re_c=10^5$, small enough that the boundary layer is likely laminar upstream of
the separation point. The average separation bubble length is
$L_{\text{sep}}=0.2c$. More information regarding the separation system and the flat plate model can be found in~\cite{DeemAIAA2017}.

PIV velocity data is sampled at $f_s=1600$\,Hz, with a resolution of $319 \times 62$ pixels. 
 The PIV vorticity dataset for the separated flow studied here consists of $m = 3000$ snapshot pairs (the training data), with a state dimension $n= 319 \times 62=19778$.  We also take another $m_{\text{test}}=3000$ snapshot pairs as testing data.

 \begin{figure}[thb]
        \centering
        \begin{subfigure}[b]{0.45\textwidth}
                    \centering
    \begin{tikzonimage}[width=.9\linewidth]{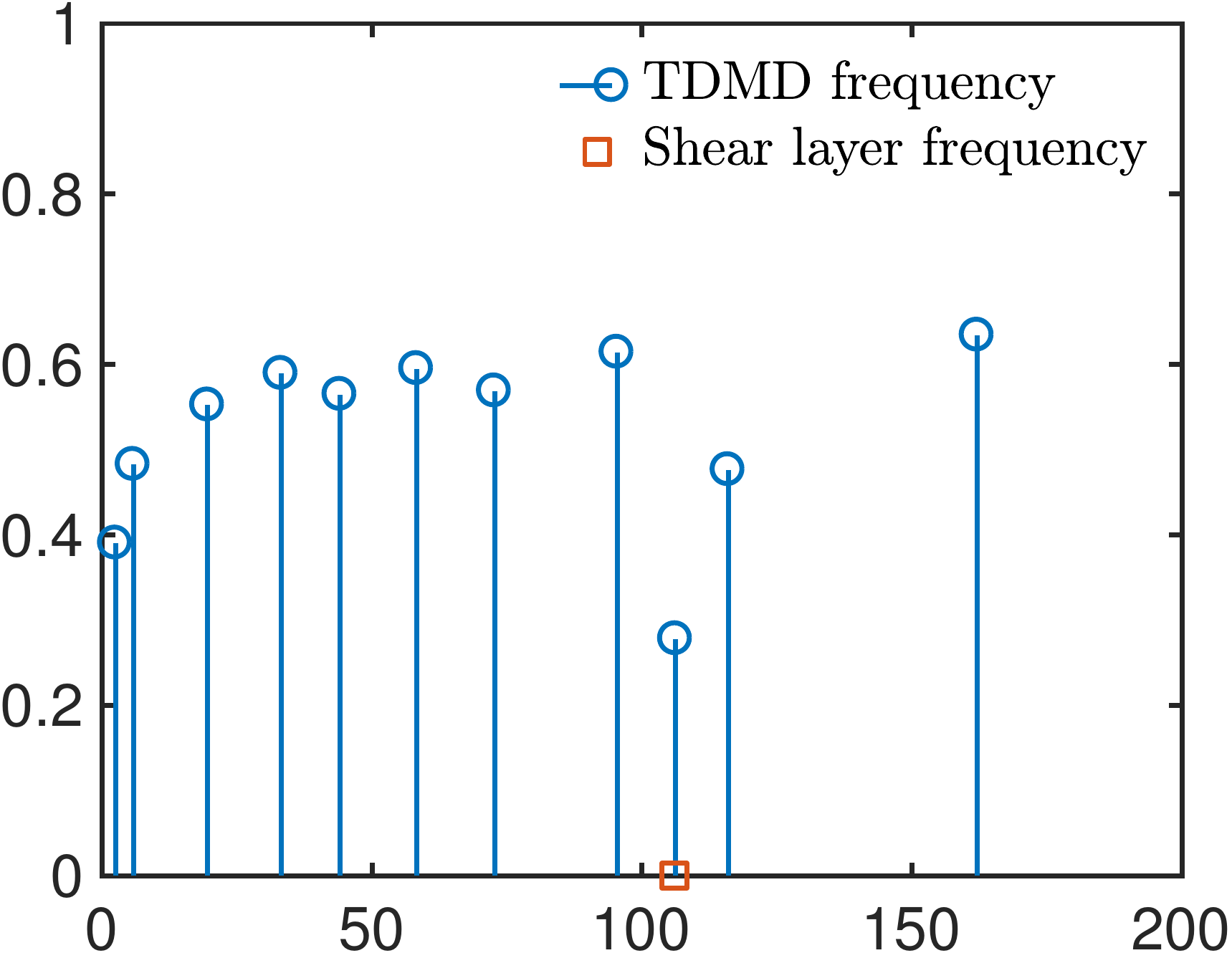}
      \small
      \node[below] at (0.5,0) {Frequency (Hz)};
      \node[anchor=south,rotate=90] at (0,0.5) {$\alpha$};
    \end{tikzonimage}
        \caption{TDMD, accuracy criterion}
        \end{subfigure}%
        \begin{subfigure}[b]{0.46\textwidth}
                    \centering
    \begin{tikzonimage}[width=.9\linewidth]{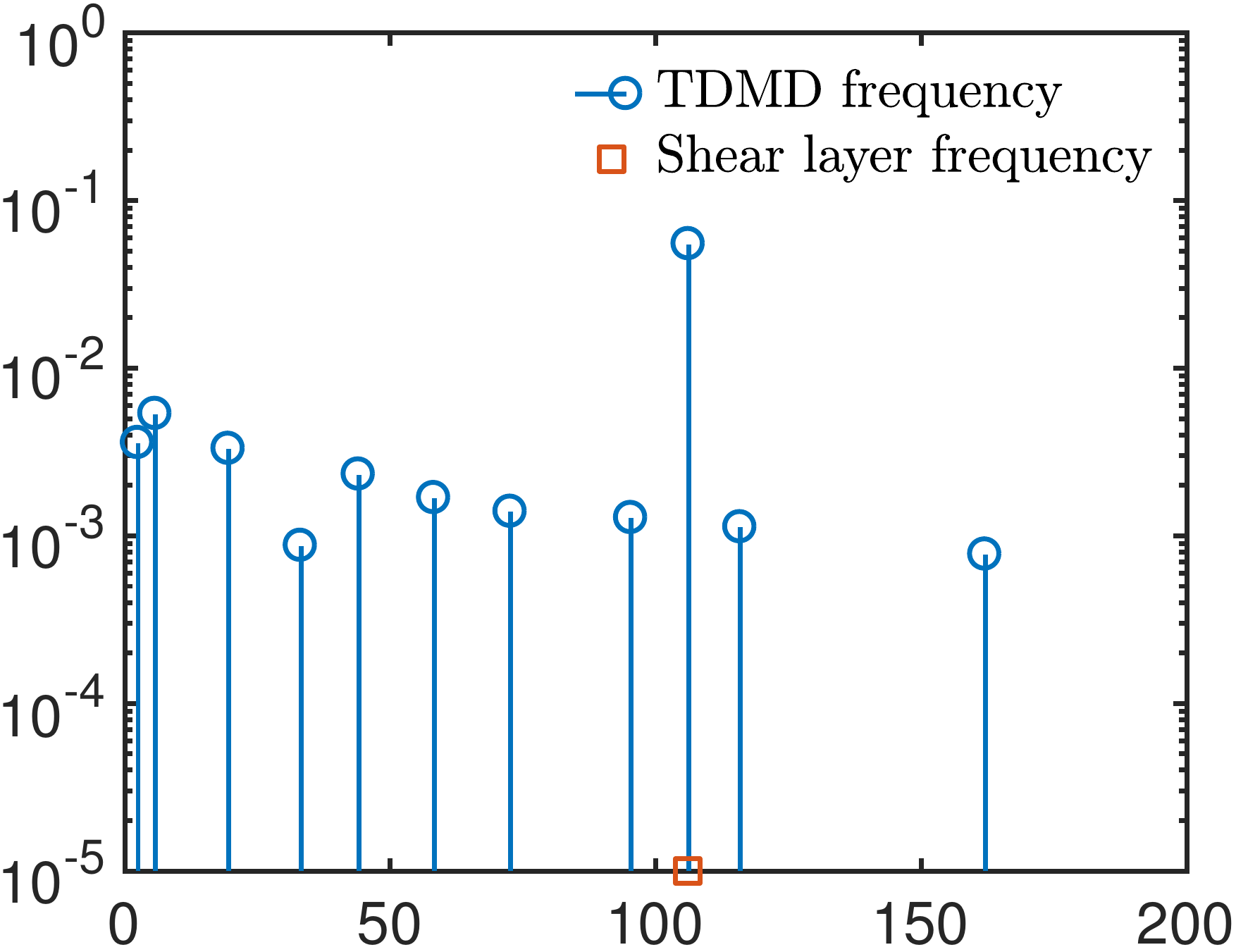}
      \small
      \node[below] at (0.5,0) {Frequency (Hz)};
      \node[anchor=south,rotate=90] at (0,0.5) {$\beta$};
    \end{tikzonimage}
                \caption{TDMD, mode amplitude}
        \end{subfigure}%
        \\
          \begin{subfigure}[b]{0.9\textwidth}
    \includegraphics[width=0.95\linewidth]{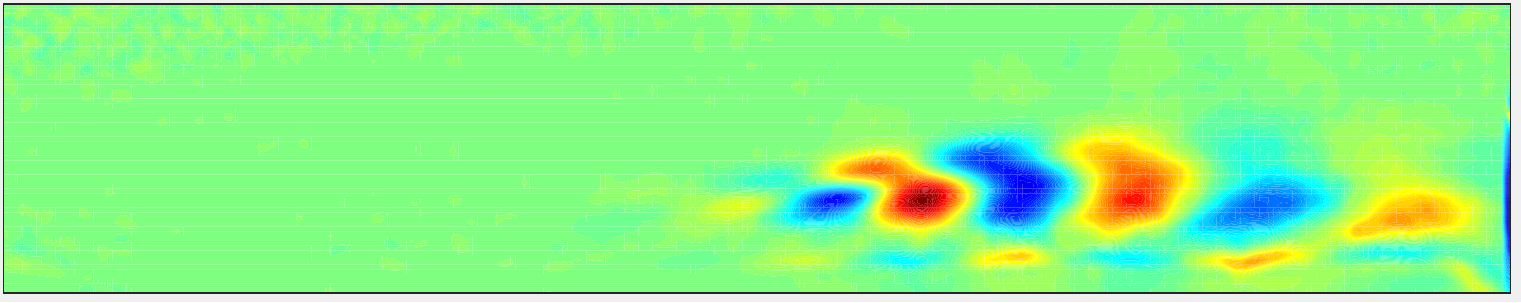}
        \centering
    \caption{TDMD mode, $f=106$ Hz}
  \end{subfigure}%
  \caption{TDMD frequency ($f_{\text{TDMD}}$) and corresponding mode
    error/amplitude. Mode amplitudes are normalized by the maximum mode amplitude. The truncation level is $r=25$. The shear layer frequency
    $f_{\text{SL}}=106$ Hz is denoted with a red square, and corresponds to the
    most accurate (smallest~$\alpha$) and largest amplitude (largest~$\beta$)
    mode.}
  \label{fig:SeparationDMD}
\end{figure}

This particular experimental dataset has been used and studied in previous work
\cite{hemati2016improving}, in which the shear layer frequency was found to be $f_{\text{SL}}=106$\,Hz. The shear layer frequency is a periodic roll-up of the shear layer due to the Kelvin-Helmholtz instability. The shear layer frequency $f_{\text{SL}}$ can be identified by applying total-least-squares DMD (TDMD), a variant of DMD which makes use of total-least-square regression to improve the accuracy of DMD for noisy data \cite{dawson2016characterizing, hemati2015biasing}.
As in~\cite{hemati2016improving}, we use a truncation level of $r=25$, which corresponds to preserving $74\%$ of the energy of the data. In this example the time spacing is $\triangle t=1/f_s=(1/1600)s$.

For comparison, we also compute the time-averaged mode amplitude $\beta$, as in
the example in section~\ref{sec:cylinder} (e.g., Figure~\ref{fig:Re413DMD}(b)).
The DMD frequencies are plotted against their accuracy criterion values and mode
amplitudes in Figure \ref{fig:SeparationDMD} (a)--(b). It is observed that
$f_{\text{SL}}=106$ Hz is accurately identified by TDMD. In addition, it stands
out by having a small mode error. The DMD mode associated with shear layer
frequency is plotted in Figure \ref{fig:SeparationDMD}(c), and it agrees with
the mode identified in previous work \cite{hemati2016improving}.

\paragraph{KDMD}
\begin{figure}[thb]
        \centering
\begin{subfigure}[b]{0.45\textwidth}
                    \centering
    \begin{tikzonimage}[width=.9\linewidth]{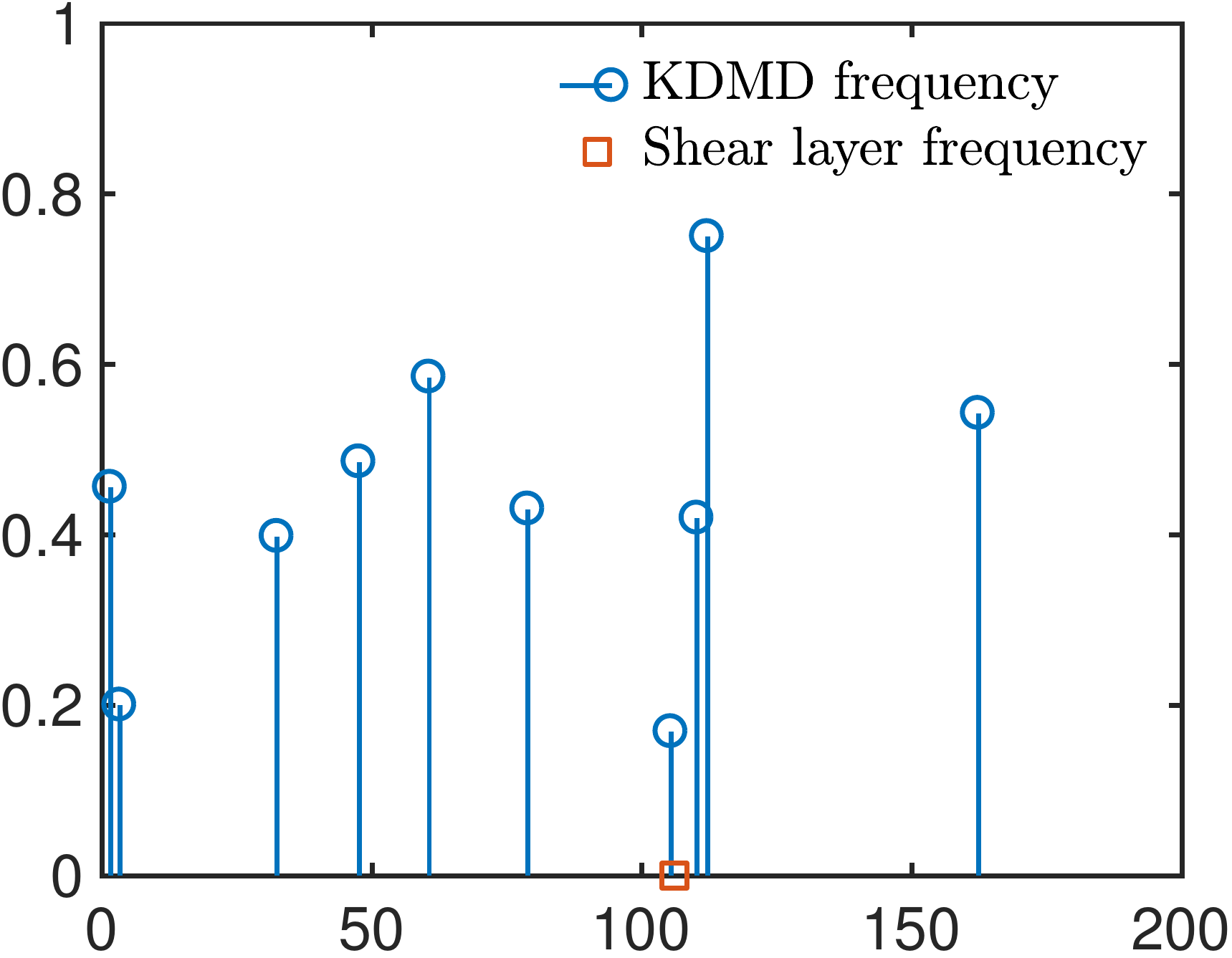}
      \small
      \node[below] at (0.5,0) {Frequency (Hz)};
      \node[anchor=south,rotate=90] at (0,0.5) {$\alpha$};
    \end{tikzonimage}
        \caption{KDMD, accuracy criterion}
        \end{subfigure}%
        \begin{subfigure}[b]{0.46\textwidth}
                    \centering
    \begin{tikzonimage}[width=.9\linewidth]{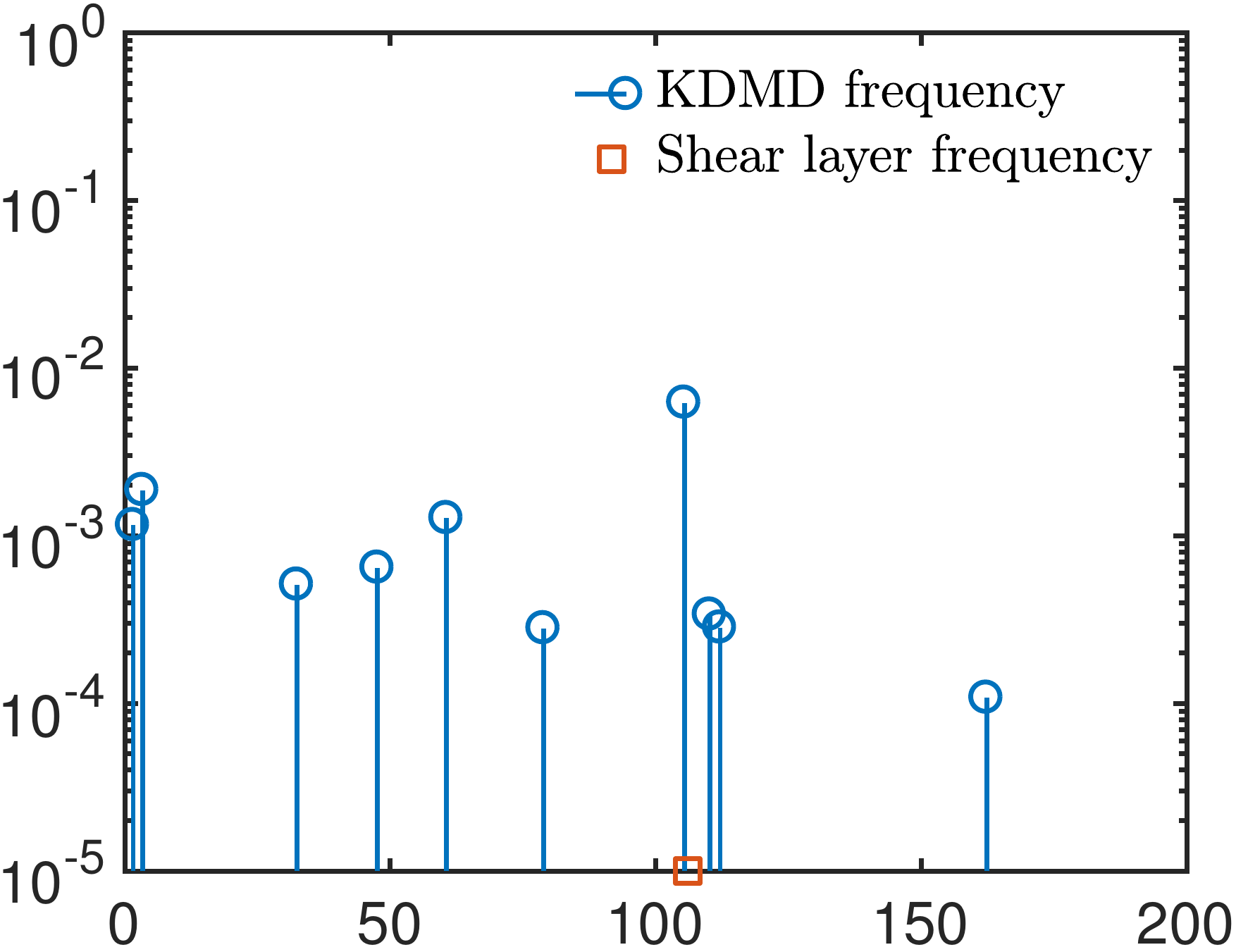}
      \small
      \node[below] at (0.5,0) {Frequency (Hz)};
      \node[anchor=south,rotate=90] at (0,0.5) {$\beta$};
    \end{tikzonimage}
                \caption{KDMD, mode amplitude}
        \end{subfigure}%
        \\
          \begin{subfigure}[b]{0.9\textwidth}
    \includegraphics[width=0.95\linewidth]{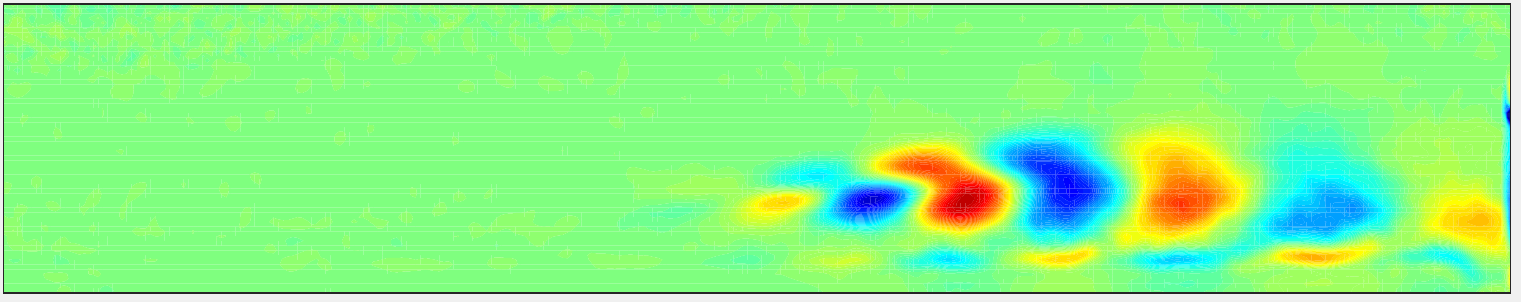}
        \centering
    \caption{KDMD mode, $f=105$ Hz}
  \end{subfigure}%
	\caption{KDMD frequency ($f_{\text{KDMD}}$) and corresponding mode error/amplitude. The truncation level is $r=25$. The shear layer frequency $f_{\text{SL}}=106$ Hz is denoted with a red square.}
	\label{fig:SeparationKDMD}
\end{figure}

We apply KDMD to this dataset, using polynomial kernels of degree $d=5$, again
with a truncation level of $r=25$. Eigenvalue frequencies, and corresponding
accuracy criterion values and mode amplitudes are plotted in
Figure~\ref{fig:SeparationKDMD}(a)--(b). We observe that the shear layer
frequency has small error and large mode amplitude, and once again verify that
the DMD mode associated with shear layer frequency
(Figure~\ref{fig:SeparationKDMD}(c)) agrees closely with that found in previous
work~\cite{hemati2016improving}.

\section{Conclusion and outlook}
Exploiting the connection between DMD and the Koopman operator, we have
presented an accuracy criterion to evaluate the quality (accuracy) of Koopman
eigenpairs approximated with DMD variants. The criterion does not assume access
to the analytical Koopman spectral decomposition, which is generally unknown in
practice. Furthermore, the proposed accuracy criterion naturally applies to
other variants of DMD, because it is based on the general notion of Koopman
eigenfunctions. The proposed accuracy criterion is validated with an synthetic
system where the analytical Koopman eigenpairs are known. Using this the
accuracy criterion, we present a study of the performance of various kernels,
and assess their sensitivity to noisy data. In our examples,  the polynomial kernel (with finite-dimensional observables) performs well both in the sense of accuracy and robustness to noise. Exponential, Gaussian, and Laplacian kernel are able to span an infinite-dimensional function space, but the tradeoff is that they are significantly more sensitive to noise in the dataset. We demonstrate that the accuracy criterion can assist in identifying accurate and physically relevant DMD modes/eigenvalues from noisy experimental data.
The accuracy criterion is conceptually simple and easy to use. As a data-driven algorithm, depending on the nature of the problem, sometimes DMD produces relevant results, and sometimes outputs numerical artifacts. For reduced order modeling based on DMD/Koopman modes, it is hence important to assess the quality of DMD results.

Note that our proposed accuracy criterion requires that some portion of data
snapshots be kept out of the DMD analysis for purposes of assessing mode
accuracy. However, it would be possible to incorporate this additional data into
the DMD analysis after the DMD modes and eigenvalues of interest have been
identified.

The demand for accurate reduced order models (ROM) has increased rapidly in
recent years, but it is still unclear how to select a subset of Koopman
eigenpairs such that the original (nonlinear) system is accurately
approximated. In order to build any meaningful ROM we need to at least assess
the accuracy and importance of DMD-approximated Koopman eigenpairs. The present
work has shed some light on the accuracy side. However, how to select the most
dynamically important Koopman eigenpairs remains an open question. Unlike
techniques such as Proper
Orthogonal Decomposition,
in which the modes are orthogonal by construction, Koopman eigenfunctions are in
general not orthogonal (though orthogonal DMD-like modes may be obtained
\cite{noack2016recursive}). Mode amplitudes obtained by a projection of data
onto DMD modes are not necessarily always a meaningful criterion for evaluating
importance, as demonstrated in the example in section~\ref{sec:crit}. It would
be desirable to develop an importance criterion that can guide the selection of
modes for the purpose of representing the dynamics accurately.

\section*{Acknowledgments}
The authors gratefully acknowledge Dr. Jessica Shang for the experimental data
of the cylinder flow. Hao Zhang thanks Dr.\ Matthew O. Williams for the generous
guidance and help as a lab mate. This material is based upon work supported by
the Air Force Office of
Scientific Research (AFOSR) under award number FA9550-14-1-0289,
and DARPA award HR0011-16-C-0116.

\bibliography{ref}

\end{document}